\documentclass[11pt]{amsart}
\usepackage{amsfonts}

\usepackage{graphicx}
\usepackage{amscd}

\newtheorem{theorem}{Theorem}[section]
\theoremstyle{plain}

\newtheorem{lemma}[theorem]{Lemma}
\newtheorem{proposition}[theorem]{Proposition}
\newtheorem{remark}{Remark}[section]

\numberwithin{equation}{section}

\begin{document}
\title[Functional KMT for partial sums ]{A functional Hungarian construction for sums of independent random variables}
\author{Ion Grama}
\address[Grama, I.]{\\
Universit\'e de Bretagne-Sud, Laboratoire SABRES \\
Rue Ives Mainguy\\
56000 Vannes, France}
\email{ion.grama@univ-ubs.fr}
\thanks{\noindent Research of the first author supported by Deutsche
Forschungsgemeinschaft, by the German Academic Exchange Service (DAAD) and
partly by the Alexander von Humboldt Schtiftung.}
\thanks{\noindent \;\;\;\; Research of the second author supported by NSF Grant
DMS-0072162 }
\author{Michael Nussbaum}
\address[Nussbaum, M.]{\\
Department of Mathematics, Malott Hall, Cornell University\\
Ithaca, NY 14853-4201, USA}
\email{nussbaum@math.cornell.edu}
\date{March 20 2001 }
\dedicatory{\textit{Universit\'e de Bretagne-Sud and Cornell University} }

\begin{abstract}
We develop a Hungarian construction for the partial sum process of
independent non-identically distributed random variables. The process is
indexed by functions $f$ from a class $\mathcal{H}$, but the supremum over $%
f\in $ $\mathcal{H}$ is taken outside the probability. This form is a
prerequisite for the Koml\'{o}s-Major-Tusn\'{a}dy inequality in the space of
bounded functionals $l^{\infty }(\mathcal{H})$, but contrary to the latter
it essentially preserves the classical $n^{-1/2}\log n$ approximation rate
over large functional classes $\mathcal{H}$ such as the H\"{o}lder ball of
smoothness $1/2$. This specific form of a strong approximation is useful for
proving asymptotic equivalence of statistical experiments.
\end{abstract}

\maketitle

\section{Introduction\label{sec:Intro}}

Let $X_{i},$ $i=1,...,n,$ be a sequence of independent random variables with
zero means and finite variances. Let $\mathcal{H}$ be a class of real valued
functions on the unit interval $[0,1]$ and $t_{i}=i/n,$ $i=1,...,n.$ The 
\emph{partial sum process indexed by functions} is the process 
\begin{equation*}
X^{n}(f)=n^{-1/2}\sum_{i=1}^{n}f(t_{i})X_{i},\;f\in \mathcal{H}.
\end{equation*}
Suppose $f\in \mathcal{H}$ are uniformly bounded; then $X^{n}=\{X^{n}(f),$ $%
f\in \mathcal{H}\}$ may be regarded as a random element with values in $%
l^{\infty }(\mathcal{H})$- the space of real valued functionals on $\mathcal{%
H}.$ The class $\mathcal{H}$ is Donsker if $X^{n}$ converges weakly in $%
l^{\infty }(\mathcal{H})$ to a Gaussian process. We are interested in
associated coupling results, i.e.\ in finding versions of $X^{n}$ and of
this Gaussian process on a common probability space which are close as
random variables. The standard coupling results of the type ''nearby
variables with nearby laws'' (cf.\ Dudley \cite{Dud}, Section 11.6)
naturally refer to the $\sup $-metric in $l^{\infty }(\mathcal{H}):$ for an
appropriate version of $X^{n}$ ($\widetilde{X}^{n}=\{\widetilde{X}^{n}(f),$ $%
f\in \mathcal{H}\},$ say) and of a Gaussian process $\widetilde{N}^{n}=\{%
\widetilde{N}^{n}(f),$ $f\in \mathcal{H}\},$ we have 
\begin{equation}
P^{\ast }\left( \sup_{f\in \mathcal{H}}\left| \widetilde{X}^{n}\left(
f\right) -\widetilde{N}^{n}\left( f\right) \right| >x\right) \rightarrow
0,\quad x>0,  \label{In.01}
\end{equation}
where $P^{\ast }$ is the outer probability on the common probability space
(cf.\ van der Vaart and Wellner \cite{VW}, 1.9.3, 1.10.4). Here we shall
consider a different type of coupling. We are looking for versions $%
\widetilde{X}^{n},$ $\widetilde{N}^{n}$ such that 
\begin{equation}
\sup_{f\in \mathcal{H}}P\left( \left| \widetilde{X}^{n}(f)-\widetilde{N}%
^{n}(f)\right| >x\right) \rightarrow 0,\quad x>0,  \label{In.02}
\end{equation}
and such that additional exponential bounds of the
Koml\'{o}s-Major-Tusn\'{a}dy type are valid. Note that (\ref{In.02}) is
weaker than (\ref{In.01}) since the supremum is taken outside the
probability. More specifically we are interested in a construction involving
also a rate sequence $r_{n}\rightarrow 0$ such that 
\begin{equation}
\sup_{f\in \mathcal{H}}P\left( r_{n}^{-1}\left| \widetilde{X}^{n}(f)-%
\widetilde{N}^{n}(f)\right| >x\right) \leq c_{0}\exp \{-c_{1}x\},\quad x>0.
\label{In.1}
\end{equation}
Here $c_{0},$ $c_{1}$ are constants depending on the class $\mathcal{H}.$

The classical results of Koml\'{o}s, Major and Tusn\'{a}dy \cite{KMT1} and 
\cite{KMT2} refer to a $\sup $ inside the probability for $\mathcal{H}=%
\mathcal{H}_{0},$ where $\mathcal{H}_{0}$ is the class of indicators $f(t)=%
\mathbf{1}(t\leq s),$ $s\in \lbrack 0,1].$ The following bound was
established: for $r_{n}=n^{-1/2}$%
\begin{equation}
P\left( r_{n}^{-1}\sup_{f\in \mathcal{H}_{0}}\left| \widetilde{X}^{n}\left(
f\right) -\widetilde{N}^{n}\left( f\right) \right| >x\right) \leq c_{0}\exp
\{-c_{1}x\},\quad x\geq c_{2}\log n,  \label{In.1a}
\end{equation}
provided $X_{1},...,X_{n}$ is a sequence of i.i.d.\ r.v.'s fulfilling
Cram\'{e}r's condition 
\begin{equation}
E\exp \{tX_{i}\}<\infty ,\quad \left| t\right| \leq t_{0},\quad i=1,...,n,
\label{I.3}
\end{equation}
where $c_{0},$ $c_{1},$ $c_{2}$ are constants depending on the common
distribution of the\ $X_{i}$. Note that $r_{n}$ in (\ref{In.1a}) can be
interpreted as a rate of convergence in the CLT over $l^{\infty }(\mathcal{H}%
_{0})$. \emph{The main reason for a construction with the supremum outside
the probability is that an extension of (\ref{In.1a}) to larger functional
classes }$\mathcal{H}$\emph{\ in general implies a substantial loss of
approximation rate }$r_{n}$ (cp. Koltchinskii (\cite{Kol}, theorem 11.1).
Our goal is a construction where the almost $n^{-1/2}$-rate of the original
KMT result is preserved despite the passage to large functional classes\emph{%
\ }$\mathcal{H}$ like Lipschitz classes.

Couplings of the type (\ref{In.1}) have first been obtained by Koltchinskii
( \cite{Kol}, theorem 3.5) and Rio \cite{Rio94} for the empirical process of
i.i.d.\ random variables, as intermediate results. They can be extended to a
full functional KMT result, i.e.\ to a coupling in $l^{\infty }(\mathcal{H})$
with exponential bounds, but an additional control of the size of the
functional class $\mathcal{H}$ is required, usually in terms of entropy
conditions. A reduced approximation rate\emph{\ }$r_{n}$ may occur as a
result.

We carry over the functional strong approximation result from the empirical
process to the partial sum process under very general conditions: the
distributions of $X_{i}$ are allowed to be \emph{nonidentical} and \emph{%
nonsmooth}. That setting substantially complicates the task of a Hungarian
construction. We can rely on the powerful methodology of Sakhanenko \cite
{Sakh}, who established the classical coupling (\ref{In.1a}) for
nonidentical and nonsmooth summands. We stress however that for the
functional version (\ref{In.1}) we need to perform the construction entirely
anew. Our results relate to Sakhanenko's \cite{Sakh} as Koltchinskii's
theorem 3.5 relates to Koml\'{o}s, Major and Tusn\'{a}dy \cite{KMT1} and 
\cite{KMT2}.

Further motivational discussion can be grouped under headings A)-C) below.

\textbf{A) Statistical applications. }The Koml\'{o}s-Major-Tusn\'{a}dy
approximation has recently found an application in the asymptotic theory of
statistical experiments. In \cite{Nuss93} the classical KMT inequality for
the empirical process was used to establish that a nonparametric experiment
of i.i.d.\ observation on an interval can be approximated, in the sense of
Le Cam's deficiency distance, by a sequence of signal estimation problems in
Gaussian white noise. The two sequences of experiments are then
asymptotically equivalent for all purposes of statistical decision with
bounded loss. This appears as a generalization of Le Cam's theory of local
asymptotic normality, applicable to ill-posed problems like density
estimation. In particular it implies a nonparametric version of the
H\`{a}jek-Le Cam asymptotic minimax theorem. The control of the Le Cam
distance is given by a relation to likelihood processes (see Le Cam and Yang 
\cite{LeCam-Yang}). Assume that there is an element $f_{0}\in \Sigma $ such
that the measures in the experiments $\mathcal{E}^{n}$ and $\mathcal{G}^{n}$
are absolutely continuous w.r.t. $P_{f_{0}}^{n}$ and $Q_{f_{0}}^{n}$
respectively. If there are versions $d\widetilde{P}_{f}^{n}/d\widetilde{P}%
_{f_{0}}^{n}$ and $d\widetilde{Q}_{f}^{n}/d\widetilde{Q}_{f_{0}}^{n}$ of the
likelihood ratios $dP_{f}^{n}/dP_{f_{0}}^{n}$ and $dQ_{f}^{n}/dQ_{f_{0}}^{n}$
on a common probability space $(\Omega ^{n},\mathcal{F}^{n},P^{n}),$ then 
\begin{equation*}
\Delta \left( \mathcal{E}^{n},\mathcal{G}^{n}\right) \leq \sqrt{2}\sup_{f\in
\Sigma }E_{P}^{n}\left( \sqrt{d\widetilde{P}_{f}^{n}/d\widetilde{P}%
_{f_{0}}^{n}}-\sqrt{d\widetilde{Q}_{f}^{n}/d\widetilde{Q}_{f_{0}}^{n}}%
\right) ^{2}
\end{equation*}
(here the expected value on the right side coincides with the Hellinger
distance between $\widetilde{P}_{f}^{n}$ and $\widetilde{Q}_{f}^{n}$). Thus
asymptotic equivalence of experiments $\mathcal{E}^{n}$ and $\mathcal{G}^{n}$
requires a ''good'' coupling of the corresponding likelihood ratios $%
dP_{f}^{n}/dP_{f_{0}}^{n}$ and $dQ_{f}^{n}/dQ_{f_{0}}^{n}$ on a common
probability space. This is achieved by constructing the linear terms (in $%
f-f_{0}$) in the expansions of the log-likelihoods such that they are close
as random variables; hence the demand for an inequality (\ref{In.1}) with
the supremum outside the probability.\textbf{\ }

The Hungarian construction had been applied in statistics before, mostly for
results on strong approximation of particular density and regression
estimators (cf.\ Cs\"{o}rg\H{o} and R\'{e}v\'{e}sz \cite{Cs-R}). It is
typical for these results that the\textit{\ ''}supremum inside the
probability'' is needed; for such an application of the functional KMT cf.
Rio \cite{Rio94}. However for asymptotic equivalence of experiments, it
turned out that it is sufficient, and indeed preferable, to have a coupling
like (\ref{In.1}) with the ''supremum outside the probability''. Applying
theorem 3.5 of Koltchinskii \cite{Kol}, it became possible in \cite{Nuss95}
to extend the scope of asymptotic equivalence, for the density estimation
problem, down to the limit of smoothness $1/2$. Analogously the present
result is essential for establishing asymptotic equivalence of smooth
nongaussian regression models to a sequence of Gaussian experiments, cf.
Grama and Nussbaum \cite{Gram-Nuss-2}. The original result of Koml\'{o}s,
Major and Tusn\'{a}dy on the partial sum process \cite{KMT1} can be used for
asymptotic equivalence in regression models, but presumably with a
non-optimal smoothness limit as in \cite{Nuss93}.

\textbf{B) Nonidentical and nonsmooth distributions. }The assumption of
identically distributed r.v.'s substantially restricts the scope of
application of the classical KMT inequality for partial sums. However this
assumption happens to be an essential point in the original proof by
Koml\'{o}s, Major and Tusn\'{a}dy and also in much of the subsequent work.
The original bound was extended and improved by many authors.
Multidimensional versions were proved by Einmahl \cite{Einm} and Zaitsev 
\cite{Zaits1}, \cite{Zaits2} with a supremum over the class of indicators $%
\mathcal{H}_{0}$. A transparent proof of the original result was given by
Bretagnolle and Massart \cite{Bret-Mass}. We would like to mention the
series of papers by Massart \cite{Massart} and Rio \cite{Rio93a}, \cite
{Rio93b}. They treat the case of $\mathbf{R}^{k}$-valued r.v.'s $X_{i},$
indexed in $\mathbb{Z}_{+}^{d}$ with a supremum taken over classes $\mathcal{%
H}$ of indicator functions $f=\mathbf{1}_{S}$ of Borel sets $S$ satisfying
some regularity conditions. Condition (\ref{I.3}) is also relaxed to moment
assumptions, but identical distributions are still assumed.

Although there are no formal restrictions on the distributions of $X_{i}$
when performing a Hungarian construction, it is not possible to get the
required closeness between the constructed r.v.'s $\widetilde{X}_{i}\overset{%
d}{=}X_{i}$ and their normal counterparts $N_{i}$ if the r.v.'s $X_{i}$ are
non-identically and non-smoothly distributed (see section \ref{sec:QT}.)
This can be argued in the following way (see Sakhanenko \cite{Sakh}). Let us
consider the sum $S=X_{1}+...+X_{n},$ where $X_{i}$ takes values $\pm
(1+2^{-i}).$ Then we can identify each realization $X_{i}$ by knowing only $%
S.$ In the dyadic Hungarian scheme, the conditional distribution of $%
X_{1}+...+X_{[n/2]}$ given $S$ is considered and used for coupling with a
Gaussian random variable. However this distribution is now degenerate and
hence not useful for coupling. This problem does not appear in the i.i.d.\
case, due to the exchangeability of the\ $X_{i}$.

We adopt a method to overcome this difficulty proposed by Sakhanenko \cite
{Sakh}. In his original paper Sakhanenko treats the case of independent
non-identically distributed r.v.'s for a class of indicators of intervals $%
\mathcal{H}=\mathcal{H}_{0}$. Here we consider the problem in another
setting: $\mathcal{H}=\mathcal{H}(1/2,L)$ where $\mathcal{H}(1/2,L)$ is a
H\"{o}lder ball with exponent $1/2$ and the $\sup $ is outside the
probability, i.e.\ we give an exponential bound for the quantity (\ref{In.1}%
) uniformly in $f$ over the set of functions $\mathcal{H}(1/2,L).$ One
complication which then appears is that the pairs $(\widetilde{X}_{i},%
\widetilde{W}_{i}),$ $i=1,...,n,$ of r.v.'s $\widetilde{X}_{i}\overset{d}{=}%
X_{i}$ and $\widetilde{W}_{i}\overset{d}{=}W_{i},$ $i=1,...,n,$ constructed
on the same probability space by the KMT method are no longer independent,
even though $\widetilde{X}_{i},$ $i=1,...,n,$ and $\widetilde{W}_{i},$ $%
i=1,...,n$ are sequences of independent r.v.'s. To deal with this we have to
develop additional properties of the Hungarian construction which are not
used in the classical setting (see Lemma \ref{L B2} for details).

\textbf{C) Coupling from marginals. }A weaker coupling of $\widetilde{X}^{n}$
and $\widetilde{N}^{n}$ can be obtained as follows. Assume for a moment that
the r.v.'s $X_{i}$ are uniformly bounded: $\left| X_{i}\right| \leq L$, $%
i=1,\ldots ,n,$ and also that $\left\| f\right\| _{\infty }\leq L,$ $f\in 
\mathcal{H}$. Take a finite collection of functions $\mathcal{H}%
_{00}=(f_{j})_{j=1,\ldots ,d}\subset \mathcal{H}$ and consider $%
Z_{i}=(f(t_{i})X_{i})_{f\in \mathcal{H}_{00}}$ as random vectors in $\mathbf{%
R}^{d}$. Reasoning as in Fact 2.2 of Einmahl and Mason \cite{Einm-Mas}
(using the result of Zaitsev \cite{Zaits3} on the Prokhorov distance between
the law of $\sum_{i=1}^{n}Z_{i}$ and a Gaussian law) we infer that for all
such $\mathcal{H}_{00}$ there are versions $\widetilde{X}^{n}\left( f\right)
,$ $\widetilde{N}^{n}\left( f\right) ,$ $f\in \mathcal{H}_{00}$ (depending
on $x$) such that 
\begin{equation}
P\left( n^{1/2}\max_{f\in \mathcal{H}_{00}}\left| \widetilde{X}^{n}(f)-%
\widetilde{N}^{n}(f)\right| \geq x\right) \leq c_{0}\exp
(-c_{1}xL^{-2}),\quad x\geq 0.  \label{coupmarg}
\end{equation}
This yields (\ref{In.1}) with rate $r_{n}=n^{-1/2}$ for every finite class $%
\mathcal{H}_{00}\subset \mathcal{H}$ of size $d,$ but with constants $%
c_{0},c_{1}$ depending on $d$. Hence any attempt to construct $\widetilde{X}%
^{n}\left( f\right) $ and $\widetilde{N}^{n}\left( f\right) ,$ on the full
class $\mathcal{H}$ from (\ref{coupmarg}) is bound to entail a substantial
loss in rate $r_{n}$; but laws of the iterated logarithm can be established
in this way (cf.\ Einmahl and Mason \cite{Einm-Mas}). Thus, to obtain (\ref
{In.1}) for $r_{n}=n^{-1/2}\log ^{2}n$ and a full H\"{o}lder class $\mathcal{%
H}(1/2,L)$, the shortcut via (\ref{coupmarg}) appears not feasible, and we
revert to a direct KMT-type construction.

In order to keep the proof somewhat transparent we do not look for optimal
logarithmic terms, but we believe that the optimal rate can be obtained by
using the very delicate technique of the paper \cite{Sakh}. The main idea
is, roughly speaking, to consider some \textit{smoothed} sequences of r.v.'s
instead of the initial \textit{unsmoothed} sequence $X_{1},...,X_{n},$ and
to apply the KMT construction for the smoothed sequences. This we perform by
substituting normal r.v.'s $N_{i}$ for the original r.v.'s $X_{i},$ for even
indices $i=2k$ in the initial sequence. Thus we are able to construct one
half of our sequence and combine it with a Haar expansion of the function $%
f. $ For the other half we apply the same argument which leads to a
recursive procedure. It turns out that this kind of smoothing is enough to
obtain ''good'' quantile inequalities although it gives rise to an
additional $\log n$ term. On the other hand the usual smoothing technique
(of each r.v.\ $X_{i}$ individually) fails. Unfortunately even the above
smoothing procedure applied with normal r.v.'s is not sufficient to obtain
the best power for the $\log n$ in the KMT inequality for non-identically
distributed r.v.'s. An optimal approach is developed in the paper of
Sakhanenko \cite{Sakh} and uses r.v.'s constructed in a special way instead
of the normal r.v.'s. Roughly speaking it corresponds to taking into
consideration the higher terms in an asymptotic expansion for the
probabilities of large deviations, which dramatically complicates the
problem. For more details we refer the reader to this beautiful paper.

Nevertheless we would like to point out that the additional $\log n$ term
which appears in our KMT result does not affect the eventual applications
that we have in mind, i.e. asymptotic equivalence of sequences of
nonparametric statistical experiments. We also believe that a stronger
version of this result (with a supremum inside the probability) might be of
use for constructing efficient kernel estimators in nonparametric models.
But such an extension is beyond of the scope of the paper.

%%%%%%%%%%%%%%%%%%%%%%%%%%%%%%%%%%%%%%%%%%%%%%%%%%%%%%%%%%%%%%%%%%%%%%%%%%%%%

\section{Notation and main results\label{sec:NR}}

Let $n\in \{1,2,... \}.$ Suppose that on the probability space $(\Omega
^{\prime },\mathcal{F}^{\prime },P^{\prime })$ we are given a sequence of
independent r.v.'s $X_1,...,X_n$ such that 
\begin{equation*}
E^{\prime }X_i=0,\quad C_{\min }\leq E^{\prime }X_i^2\leq C_{\max },\quad
i=1,...,n,
\end{equation*}
where $C_{\min }<C_{\max }$ are some positive absolute constants. Hereafter $%
E^{\prime }$ is the expectation under the measure $P^{\prime }.$ Assume also
that the following extension of a condition due to Sakhanenko \cite{Sakh}
holds true: 
\begin{equation}
\lambda _nE^{\prime }\left| X_i\right| ^3\exp \{\lambda _n\left| X_i\right|
\}\leq E^{\prime }X_i^2,\quad i=1,...,n,  \label{EXT-SAKH-COND}
\end{equation}
where $\lambda _n$ is a sequence of real numbers satisfying $0<\lambda
_n<\lambda ,$ $n\geq 1,$ for some positive absolute constant $\lambda .$
Along with this, assume that on another probability space $(\Omega ,\mathcal{%
F},P)$ we are given a sequence of independent normal r.v.'s $N_1,...,N_n$
such that 
\begin{equation*}
EN_i=0,\quad EN_i^2=E^{\prime }X_i^2,
\end{equation*}
for all $i=1,...,n.$ Hereafter $E$ is the expectation under the measure $P.$

Let $\mathcal{H}(1/2,L)$ be the H\"{o}lder ball with exponent $1/2,$ i.e.\
the set of real valued functions $f$ defined on the unit interval $[0,1]$
and satisfying the following conditions 
\begin{equation*}
\left| f(x)-f(y)\right| \leq L\left| x-y\right| ^{1/2},\quad \left\|
f\right\| _{\infty }\leq L/2,
\end{equation*}
where $L$ is a positive absolute constant.

Let $t_{i}=i/n,$ $i=1,...,n$ be a uniform grid in the unit interval $[0,1].$
The notation $Y\overset{d}{=}X$ for random variables means equality in
distribution. The symbol $c$ (with possible indices) denotes a generic
positive absolute constant (more precisely this means that it is a function
only of the absolute constants introduced before).

The main result of the paper is the following.

\begin{theorem}
\label{T1} Let $n\geq 2.$ A sequence of independent r.v.'s $\widetilde{X}%
_{1},...,\widetilde{X}_{n}$ can be constructed on the probability space $%
(\Omega ,\mathcal{F},P)$ such that $\widetilde{X}_{i}\overset{d}{=}X_{i},$ $%
i=1,...,n$ and 
\begin{equation*}
\sup_{f\in \mathcal{H}(1/2,L)}P\left( \left| \sum_{i=1}^{n}f(t_{i})(%
\widetilde{X}_{i}-N_{i})\right| >x\frac{\log ^{2}n}{\lambda _{n}}\right)
\leq c_{1}\exp \{-c_{2}x\},\quad x\geq 0.
\end{equation*}
\end{theorem}

\begin{remark}
In the above theorem $X_{i},$ $i=1,...,n$ are not supposed to be identically
distributed nor to have smooth distributions, although the result is new
even in the case of i.i.d.\ r.v.'s. The r.v.'s $\widetilde{X}_{1},...,%
\widetilde{X}_{n}$ constructed are functions of the r.v.'s $N_{1},...,N_{n}$
only, so that no assumptions on the probability space $(\Omega ,\mathcal{F}%
,P)$ are required other than existence of $N_{1},...,N_{n}$.
\end{remark}

\begin{remark}
The use of condition (\ref{EXT-SAKH-COND}) instead of a more familiar Cram%
\'{e}r type condition is motivated by the desire to cover also the case of
non-identically distributed r.v.'s with subexponential moments, which
corresponds to $\lambda _{n}\rightarrow 0.$ This case cannot be treated
under Cram\'{e}r's condition, but it is important since it essentially
includes the case of non-identically distributed r.v.'s with finite moments.
\end{remark}

Theorem \ref{T1} can be formulated in the following equivalent form.

\begin{theorem}
\label{T3} Let $n\geq 2.$ A sequence of independent r.v.'s $\widetilde{X}%
_{1},...,\widetilde{X}_{n}$ can be constructed on the probability space $%
(\Omega ,\mathcal{F},P)$ such that $\widetilde{X}_{i}\overset{d}{=}X_{i},$ $%
i=1,...,n,$ and for any $t$ satisfying $\left| t\right| \leq c_{1}$%
\begin{equation*}
\sup_{f\in \mathcal{H}(1/2,L)}E\exp \left\{ t\frac{\lambda _{n}}{\log ^{2}n}%
\sum_{i=1}^{n}f(t_{i})(\widetilde{X}_{i}-N_{i})\right\} \leq \exp \left\{
c_{2}t^{2}\right\} .
\end{equation*}
\end{theorem}

Let us formulate yet another equivalent version of Theorem \ref{T1}. Assume
that on the probability space $(\Omega ^{\prime },\mathcal{F}^{\prime
},P^{\prime })$ we are given a sequence of independent r.v.'s $%
X_{1},...,X_{n}$ such that for all $i=1,...,n$

\begin{equation}
E^{\prime }X_{i}=0,\quad \lambda _{n}^{2}C_{\min }\leq E^{\prime
}X_{i}^{2}\leq C_{\max }\lambda _{n}^{2},  \label{NR.1}
\end{equation}
where $C_{\min }<C_{\max }$ are positive absolute constants and $\lambda
_{n} $ is a sequence of real numbers $0<\lambda _{n}\leq 1,$ $n\geq 1.$
Assume also that the following condition due to Sakhanenko \cite{Sakh} holds
true: 
\begin{equation}
\lambda E^{\prime }\left| X_{i}\right| ^{3}\exp \{\lambda \left|
X_{i}\right| \}\leq E^{\prime }X_{i}^{2},\quad i=1,...,n,  \label{NR.2}
\end{equation}
where $\lambda $ is a positive absolute constant. Suppose that on another
probability space $(\Omega ,\mathcal{F},P)$ we are given a sequence of
independent normal r.v.'s $N_{1},...,N_{n}$ such that for $i=1,...,n$%
\begin{equation}
EN_{i}=0,\quad EN_{i}^{2}=E^{\prime }X_{i}^{2}.  \label{NR.4}
\end{equation}

\begin{theorem}
\label{T3-2} Let $n\geq 2.$ A sequence of independent r.v.'s $\widetilde{X}%
_{1},...,\widetilde{X}_{n}$ can be constructed on the probability space $%
(\Omega ,\mathcal{F},P)$ such that $\widetilde{X}_{i}\overset{d}{=}X_{i},$ $%
i=1,...,n,$ and for any $t$ satisfying $\left| t\right| \leq c_{1}$%
\begin{equation*}
\sup_{f\in \mathcal{H}(1/2,L)}E\exp \left\{ \frac{t}{\log ^{2}n}%
\sum_{i=1}^{n}f(t_{i})(\widetilde{X}_{i}-N_{i})\right\} \leq \exp \left\{
c_{2}t^{2}\right\} .
\end{equation*}
\end{theorem}

We shall give a proof of Theorem \ref{T3-2} in Section \ref{sec:PMR}.

Now we turn to a particular case of the above results. Assume that the
sequence of independent r.v.'s $X_{1},...,X_{n}$ is such that 
\begin{equation}
E^{\prime }X_{i}=0,\quad C_{\min }\leq E^{\prime }X_{i}^{2}\leq C_{\max
},\quad i=1,...,n,  \label{NR.7}
\end{equation}
for some positive absolute constants $C_{\min }<C_{\max }.$ Assume also that
the following Cram\'{e}r type condition holds true: 
\begin{equation}
E^{\prime }\exp \{C_{1}\left| X_{i}\right| \}\leq C_{2},\quad i=1,...,n,
\label{NR.8}
\end{equation}
where $C_{1}$ and $C_{2}$ are positive absolute constants.

\begin{theorem}
Let $n\geq 2.$ A sequence of independent r.v.'s $\widetilde{X}_{1},...,%
\widetilde{X}_{n}$ can be constructed on the probability space $(\Omega ,%
\mathcal{F},P)$ such that $\widetilde{X}_{i}\overset{d}{=}X_{i},$ $i=1,...,n$
and 
\begin{equation*}
\sup_{f\in \mathcal{H}(1/2,L)}P\left( \left| \sum_{i=1}^{n}f(t_{i})(%
\widetilde{X}_{i}-N_{i})\right| >x\log ^{2}n\right) \leq c_{1}\exp
\{-c_{2}x\},\quad x\geq 0.
\end{equation*}
\end{theorem}

To deduce this result from Theorem \ref{T1}, it suffices to note that
Sakhanenko's condition (\ref{NR.2}) holds true with $\lambda _{n}=const$
depending on $C_{\min },$ $C_{1}$ and $C_{2},$ under (\ref{NR.7}) and (\ref
{NR.8}).

\begin{remark}
\label{Rem-S-for-norm}It should be mentioned that Sakhanenko's condition (%
\ref{NR.2}) holds true for the normal r.v.'s $N_{1},...,N_{n}$ only if the
constant $\lambda $ is small enough, namely if $\lambda \leq c\left(
EN_{i}^{2}\right) ^{-1/2}.$ Since the function $\alpha |x|^{3}\exp (\alpha
|x|)$ is increasing in $\alpha ,$ the condition (\ref{NR.2}) holds true for
any $\lambda \leq \lambda ^{\prime }$ if it holds true with some $\lambda
=\lambda ^{\prime }$. Therefore without loss of generality it can be assumed
that the constant $\lambda $ fulfills $\lambda \leq c/C_{\max }\leq c\left(
E^{\prime }X_{i}^{2}\right) ^{-1/2},$ $i=1,...,n,$ thus ensuring that (\ref
{NR.2}) holds true also for $N_{1},...,N_{n}.$
\end{remark}

\section{Elementary properties of Haar expansions}

For the following basic facts we refer to Kashin and Saakyan \cite{Kash-Saak}%
). The Fourier-Haar basis on the interval $[0,1]$ is introduced as follows.
Consider the dyadic system of partitions by setting 
\begin{equation*}
s_{k,j}=j2^{-k},
\end{equation*}
for $\,j=1,...,2^{k}\,$ and 
\begin{equation}
\Delta _{k,1}=[0,s_{k,1}],\quad \Delta _{k,j}=(s_{k,j-1},s_{k,j}],
\label{HE.2x}
\end{equation}
for $\,j=2,...,2^{k},$ where $\,k\geq 0.$ Define Haar functions via
indicators $1(\Delta _{k,j})$ 
\begin{equation}
h_{0}=1(\Delta _{0,1}),\quad h_{k,j}=2^{k/2}(1(\Delta _{k+1,2j-1})-1(\Delta
_{k+1,2j})),  \label{HE.2x-1}
\end{equation}
for $j=1,...,2^{k}$ and $k\geq 0.$

If $f$ is a function from $\mathcal{L}_{2}([0,1])$ then the following Haar
expansion 
\begin{equation*}
f=c_{0}(f)h_{0}+\sum_{k=0}^{\infty }\sum_{j=1}^{2^{k}}c_{k,j}(f)h_{k,j},
\end{equation*}
holds true with Fourier-Haar coefficients 
\begin{equation}
c_{0}(f)=\int_{0}^{1}f(u)h_{0}(u)du,\quad
c_{k,j}(f)=\int_{0}^{1}f(u)h_{k,j}(u)\,du,  \label{HE.4}
\end{equation}
for $j=1,...,2^{k}$ and $k\geq 0.$ Along with this, consider the truncated
Haar expansion 
\begin{equation}
f_{_{m}}=c_{0}(f)h_{0}+\sum_{k=0}^{m-1}\sum_{j=1}^{2^{k}}c_{k,j}(f)h_{k,j},
\label{HE.5}
\end{equation}
for some $m\geq 1$.

\begin{proposition}
\label{P HE1}For $f\in \mathcal{H}(1/2,L)$ we have 
\begin{equation*}
|c_{0}(f)|\leq L/2,\quad |c_{k,j}(f)|\leq 2^{-3/2}L2^{-k},
\end{equation*}
for $k=0,1,...$ and $j=1,...,2^{k}.$
\end{proposition}

\begin{proof}
It is easy to see that 
\begin{eqnarray*}
c_{k,j}(f) &=&2^{k/2}(\int_{\Delta _{k+1,2j-1}}f(u)\,du-\int_{\Delta
_{k+1,2j}}f(u)\,du), \\
&=&2^{k/2}\int_{\Delta _{k+1,2j-1}}(f(u)-f(u+2^{-(k+1)}))\,du.
\end{eqnarray*}
Since $f$ is in the H\"{o}lder ball $\mathcal{H}(\frac{1}{2},L)$ we get 
\begin{eqnarray*}
|c_{k,j}(f)| &\leq &2^{k/2}\sup_{u\in \Delta
_{k+1,2j-1}}|f(u)-f(u+2^{-(k+1)})|\int_{\Delta _{k+1,2j-1}}du \\
&\leq &2^{k/2}L2^{-(k+1)/2}2^{-(k+1)}\leq 2^{-3/2}L2^{-k}.
\end{eqnarray*}
\end{proof}

Next we give an estimate for the uniform distance between $f$ and $f_{m}.$

\begin{proposition}
\label{P HE2}For $f\in \mathcal{H}(1/2,L)$ we have 
\begin{equation*}
\sup_{0\leq t\leq 1}|f(t)-f_{m}(t)|\leq L2^{-m/2}.
\end{equation*}
\end{proposition}

\begin{proof}
It is easy to check (see for instance Kashin and Saakyan \cite{Kash-Saak},
p. 81) that, whenever $t\in \Delta _{m,j},$ 
\begin{equation*}
f_{m}(t)=2^{m}\int_{\Delta _{m,j}}f(s)\,d\,s,
\end{equation*}
for $j=1,...,2^{m}$, which gives us $f_{m}(t)=f(\widetilde{t}_{m,j})$, with
some $\widetilde{t}_{m,j}\in \Delta _{m,j}.$ Since $f(t)$ is in the H\"{o}%
lder ball $\mathcal{H}(\frac{1}{2},L)$, we obtain for any $j=1,...,2^{m}$
and $t\in \Delta _{m,j}$%
\begin{equation*}
|f(t)-f_{m}(t)|=|f(t)-f_{m}(\widetilde{t}_{m,j})|\leq L|t-\widetilde{t}%
_{m,j}|^{1/2}\leq L2^{-m/2}.
\end{equation*}
\end{proof}

\section{Background on quantile transforms \label{sec:QT}}

Let $(\Omega ^{\prime },\mathcal{F}^{\prime },P^{\prime })$ be a probability
space. Let $\lambda $ be a real number such that $0<\lambda <\infty .$
Denote by $\mathfrak{D}(\lambda )$ the set of all r.v.'s $S$ on the
probability space $(\Omega ^{\prime },\mathcal{F}^{\prime },P^{\prime })$
which can be represented as a sum $S=X_{1}+...+X_{n}$ of some independent
r.v.'s on $(\Omega ^{\prime },\mathcal{F}^{\prime },P^{\prime })$ for some $%
n\geq 1,$ satisfying relations (\ref{Q.1}), (\ref{Q.2}) below:

\begin{itemize}
\item  The r.v.'s $X_{1},...,X_{n}$ have zero means and finite variances: 
\begin{equation}
E^{\prime }X_{i}=0,\quad 0<E^{\prime }X_{i}^{2}<\infty   \label{Q.1}
\end{equation}
for any $i=1,...,n.$

\item  Sakhanenko's condition 
\begin{equation}
\lambda E^{\prime }\left| X_{i}\right| ^{3}\exp \left\{ \lambda
|X_{i}|\right\} <E^{\prime }X_{i}^{2},  \label{Q.2}
\end{equation}
is satisfied for all $i=1,...,n.$
\end{itemize}

Let $\mu $ be a real number satisfying $0<\mu <\infty .$ By $\mathfrak{D}%
_0(\lambda ,\mu )$ we denote the subset of all r.v.'s $S\in \mathfrak{D}%
(\lambda )$ which additionally satisfy the following smoothness condition (%
\ref{Q.3}):

\begin{itemize}
\item  For any $0<\varepsilon <1,$ we have 
\begin{equation}
\sup_{|h|\leq \varepsilon }\int_{|t|>\varepsilon }\left| \frac{E^{\prime
}\exp \left\{ (\mathbf{i}t+h)S\right\} }{E^{\prime }\exp \left\{ hS\right\} }%
\right| dt\leq \frac{\mu }{\varepsilon E^{\prime }S^{2}},  \label{Q.3}
\end{equation}
where $\mathbf{i=}\sqrt{-1}.$
\end{itemize}

\begin{remark}
In the sequel we shall assume that $\mu $ is a positive absolute constant,
and therefore, we shall drop the dependence on $\mu $ in the notation for $%
\mathfrak{D}_{0}(\lambda ,\mu ),$ i.e.\ we write for short $\mathfrak{D}%
_{0}(\lambda )=\mathfrak{D}_{0}(\lambda ,\mu ).$
\end{remark}

We now introduce the \textit{quantile transformation} and the associated
basic inequality (see Lemma \ref{L Q1}). Assume that on probability space $%
(\Omega ^{\prime },\mathcal{F}^{\prime },P^{\prime })$ we are given an
arbitrary r.v.\ $X$ of mean zero and finite variance: $E^{\prime }X=0$ and $%
E^{\prime }X^{2}<\infty .$ Assume that on another probability space $(\Omega
,\mathcal{F},P)$ we are given a normal r.v.\ $N$ with the same mean and
variance: $EN=0$ and $EN^{2}=E^{\prime }X^{2}.$ Let $F_{X}(x)$ and $\Phi
_{N}(x)$ be the distribution functions of $X$ and $N$ respectively. Note
that the r.v.\ $U=\Phi _{N}(N)$ is distributed uniformly on $[0,1]$. Define
the r.v. $\widetilde{X}$\ to be the solution of the equation 
\begin{equation*}
F_{X}(\widetilde{X})=\Phi _{N}(N)=U.
\end{equation*}
The r.v. $\widetilde{X}$ is called a quantile transformation of \ $N.$ It is
easy to see that a solution\ $\widetilde{X}$ always exists and has
distribution function $F,$ although it need not be unique. In the case of
non-uniqueness, we choose one of the possible solutions.

The following assertion follows from the results in Sakhanenko \cite{Sakh}
(see Theorem 4, p. 10).

\begin{lemma}
\label{L Q1} Set $B^{2}=E^{\prime }X^{2}=EN^{2}.$ In addition to the above
suppose that $X\in \mathfrak{D}_{0}(\lambda ).$ Then 
\begin{equation*}
\left| \widetilde{X}-N\right| \leq \frac{c_{1}}{\lambda }\left\{ 1+\frac{%
\widetilde{X}^{2}}{B^{2}}\right\} ,
\end{equation*}
provided $|\widetilde{X}|\leq c_{2}\lambda B^{2}$ and $\lambda B\geq c_{3},$
where $c_{1},$ $c_{2}$ and $c_{3}$ are positive absolute constants.
\end{lemma}

Let us now introduce the \textit{conditional quantile transformation} and
the associated basic inequality (Lemma \ref{L Q2} below).

Assume that on the probability space $(\Omega ^{\prime },\mathcal{F}^{\prime
},P^{\prime })$ we are given two independent r.v.'s $X_{1},$ $X_{2}$ of
means zero and finite variances: $E^{\prime }X_{i}=0$ and $E^{\prime
}X_{i}^{2}<\infty ,$ for $i=1,2.$ Assume further that on another probability
space $(\Omega ,\mathcal{F},P)$ we are given two normal r.v.'s $N_{1},$ $%
N_{2}$ with the same means and variances: $EN_{i}=0$ and $%
EN_{i}^{2}=E^{\prime }X_{i}^{2},$ for $i=1,2.$ Set $X_{0}=X_{1}+X_{2}$ and $%
N_{0}=N_{1}+N_{2}.$ Denote $B_{i}=E^{\prime }X_{i}^{2},$ $\alpha
_{1}=B_{1}/B_{2},$ $\alpha _{2}=B_{2}/B_{1}.$ Suppose that we have
constructed a \ $\widetilde{X}_{0}$ having the same distribution as $X_{0},$
and which depends only on $N_{0}$ and on some random vector $W.$ Suppose
that $N_{1}$ and $N_{2}$ do not depend on $W.$ We wish to construct $X_{1}$
and $X_{2}.$ Let $F_{T_{0}|X_{0}}(x|y)$ be the conditional distribution
function of \ $T_{0}=\alpha _{2}X_{1}-\alpha _{1}X_{2}$ given $X_{0}=y$ and $%
\Phi _{V_{0}}(x)$ be the distribution function of the normal r.v.\ $%
V_{0}=\alpha _{2}N_{1}-\alpha _{1}N_{2}.$ Define\ $\widetilde{T}_{0}$ to be
the solution of the equation 
\begin{equation*}
F_{T_{0}|X_{0}}(\widetilde{T}_{0}|\widetilde{X}_{0})=\Phi _{V_{0}}(V_{0})=U.
\end{equation*}
The r.v.\ $\widetilde{T}_{0}$ is called a conditional quantile
transformation of\ $V_{0}$ given $\widetilde{X}_{0}.$

\begin{proposition}
\label{R Q1} Set $\widetilde{X}_{1}=\alpha _{0}^{-1}\left( T_{0}+\alpha _{1}%
\widetilde{X}_{0}\right) $ and $\widetilde{X}_{2}=\alpha _{0}^{-1}\left(
T_{0}-\alpha _{2}\widetilde{X}_{0}\right) .$ Then $\widetilde{X}_{1}$ and $%
\widetilde{X}_{2}$ are independent and such that $\widetilde{X}_{1}\overset{d%
}{=}X_{1},$ $\widetilde{X}_{2}\overset{d}{=}X_{2}.$ Moreover $\widetilde{X}%
_{1}$ and $\widetilde{X}_{2}$ are functions of the r.v.'s $\widetilde{X}_{0},
$ $N_{1}$ and $N_{2}$ only.
\end{proposition}

\begin{proof}
Consider \ $U=\Phi (V_{0}).$ It is clear that the distribution of $U$ is
uniform on $[0,1].$ Since $V_{0}=\alpha _{2}N_{1}-\alpha _{1}N_{2}$ and $%
N_{0}=N_{1}+N_{2}$ are normal and uncorrelated, $U$ and $N_{0}$ are
independent. Since $(N_{1},N_{2})$ does not depend on $W,$ we conclude that $%
U$ does not depend on $N_{0}$ and $W.$ But\ $\widetilde{X}_{0}$ is a
function of $N_{0}$ and $W$ only. Hence $U$ and $\widetilde{X}_{0}$ are also
independent.

Next, since the uniform r.v.\ $U$ does not depend on $\widetilde{X}_{0},$ we
easily check that the distribution of $\widetilde{T}_{0}$ given $\widetilde{X%
}_{0}=y,$ for any real $y,$ is exactly $F_{T_{0}|X_{0}}(\cdot |y).$ Taking
into account that $\widetilde{X}_{0}\overset{d}{=}X_{0},$ we conclude that
the two-dimensional distributions of the pairs $(\widetilde{T}_{0},%
\widetilde{X}_{0})$ and $(T_{0},X_{0})$ coincide. From this we obtain in
particular that $\widetilde{X}_{1}$ and $\widetilde{X}_{2}$ are independent
and that $\widetilde{X}_{1}\overset{d}{=}X_{1},$ $\widetilde{X}_{2}\overset{d%
}{=}X_{2}.$ Moreover it is obvious from the construction that $\widetilde{X}%
_{1}$ and $\widetilde{X}_{2}$ are functions of $\widetilde{X}_{0},$ $N_{1}$
and $N_{2}$ only.
\end{proof}

The following assertion follows from the results in Sakhanenko \cite{Sakh}
(see Theorem 6, p. 20).

\begin{lemma}
\label{L Q2}Set $B=B_{1}B_{2}/B_{0}.$ In addition to the above suppose that $%
X_{1},X_{2}\in \mathfrak{D}_{0}(\lambda ).$ Then 
\begin{equation*}
\left| \widetilde{T}_{0}-V_{0}\right| \leq \frac{c_{1}}{\lambda }\frac{B_{0}%
}{B}\left\{ 1+\frac{\widetilde{T}_{0}^{2}}{B^{2}}+\frac{\widetilde{X}_{0}^{2}%
}{B^{2}}\right\} ,
\end{equation*}
provided $|\widetilde{T}_{0}|\leq c_{2}\lambda B^{2},$ $|\widetilde{X}%
_{0}|\leq c_{2}\lambda B^{2}$ and $\lambda B\geq c_{3},$ where $c_{1},$ $%
c_{2}$ and $c_{3}$ are absolute constants.
\end{lemma}

\section{A construction for non-identically distributed r.v.'s. \label%
{sec:CNST}}

In this section we assume that we are given a sequence of independent r.v.'s 
$X_{i},$ $i=1,...,n,$ satisfying (\ref{NR.1}) and (\ref{NR.2}). We shall
construct a version of this sequence and an appropriate sequence of
independent normal r.v.'s $N_{i},$ $i=1,...,n,$ on the same probability
space such that these are as close as possible.\textit{\ }More precisely,
the construction is performed so that the quantile inequalities in Section 
\ref{sec:QT} are applicable. Of course the sequences which we obtain are
dependent. To assure that this dependence remains under control, we
partition the initial sequence into dyadic blocks with similar size of
variances. Some prerequisites for this are given in the next section. The
construction itself is performed in Section \ref{sec: CONSTR}.

\subsection{A dyadic blocking procedure\label{sec:SNot}}

In this section we exhibit a special partition of the initial sequence into
dyadic blocks so that the sums of the $X_{i}$ inside the blocks at any
dyadic level have approximately the same variances. This will be used for
proving quantile inequalities in Section \ref{sec:QI} and some exponential
bounds in Section \ref{sec:PMR} (see Lemma \ref{L P2} and Proposition \ref{P
M2}).

Assume that $n>n_{\min }\geq 1,$ where $n_{\min }$ is an absolute constant
whose precise value will be indicated below. Set $M=[\log _{2}(n/n_{\min
})]. $ It is clear that $M\geq 0$ and $n_{\min }2^{M}\leq n<n_{\min
}2^{M+1}. $ Let $J_{M}=\{1,...,n\}$ and define consecutively $%
J_{m}=\{i:2i\in J_{m+1}\},$ for $m=0,...,M-1.$ Alternatively, for any $%
m=0,...,M$ the set of indices $J_{m}$ can be defined as follows: 
\begin{equation*}
J_{m}=\left\{ i:1\leq i2^{M-m}\leq n\right\} .
\end{equation*}
Let $n_{m}$ denote the last element in $J_{m}$ i.e. $n_{m}=\#J_{m}.$ It is
not difficult to see that $n_{\min }\leq n_{0}\leq 2n_{\min }.$

Recall that each r.v. $X_{i}$ is attached to a design point $t_{i}=i/n,$ $%
i=1,...,n$. Set 
\begin{equation}
t_{i}^{m}=t_{i2^{M-m}},\quad X_{i}^{m}=X_{i2^{M-m}},\quad m=0,...,M,\quad
i\in J_{m}.  \label{Xmi}
\end{equation}
Our next task is to split each sequence $X_{i}^{m},$ $i\in J_{m}$ into
dyadic blocks so that the sums of $X_{i}^{m}$ over blocks at a given
resolution level $m$ have approximately the same variances. To ensure this
we shall introduce the strictly increasing function $b_{m}(t):[0,1]%
\rightarrow \lbrack 0,1],$ which is related to the variances of $X_{i}^{m}$
as follows: 
\begin{equation*}
b_{m}(t)=\int_{0}^{t}\beta _{m}(s)ds/\int_{0}^{1}\beta _{m}(s)ds,\quad t\in
(0,1],\quad b_{m}(0)=0,
\end{equation*}
where 
\begin{equation*}
\beta _{m}(s)=\left\{ 
\begin{array}{ll}
E^{\prime }(X_{i}^{m})^{2}, & \text{if }s\in (t_{i-1}^{m},t_{i}^{m}],\;i\in
J_{m}, \\ 
E^{\prime }(X_{n_{m}}^{m})^{2}, & \text{if }s\in (t_{n_{m}}^{n},1].
\end{array}
\right.
\end{equation*}
Let $a_{m}(t)$ be the inverse of $b_{m}(t)$, i.e. 
\begin{equation}
a_{m}(t)=\inf \left\{ s\in \lbrack 0,1]:b_{m}(s)>t\right\} .  \label{BCN.4a}
\end{equation}
It is easy to see that condition (\ref{NR.1}) implies that both $b_{m}(t)$
and $a_{m}(t)$ are Lipschitz functions: for any $t_{1},t_{2}\in \lbrack
0,1], $ we have 
\begin{equation}
\left| b_{m}(t_{2})-b_{m}(t_{1})\right| \leq L_{\max }\left|
t_{2}-t_{1}\right| ,\quad \left| a_{m}(t_{2})-a_{m}(t_{1})\right| \leq
L_{\max }\left| t_{2}-t_{1}\right| ,  \label{BCN.5b}
\end{equation}
where $L_{\max }=C_{\max }/C_{\min }.$ Consider the dyadic scheme of
partitions 
\begin{equation*}
\Delta _{k,j},\quad j=1,...,2^{k},\quad k=0,...M,
\end{equation*}
of the interval $[0,1]$ as defined by (\ref{HE.2x}). For any $m=0,...,M,$
denote by $I_{k,j}^{m}$ the set of those indices $i\in J_{m}$ for which $%
b_{m}(t_{i}^{m})$ falls into $\Delta _{k,j}$, i.e.\ 
\begin{equation*}
I_{k,j}^{m}=\left\{ i\in J_{m}:b_{m}(t_{i}^{m})\in \Delta _{k,j}\right\}
,\quad j=1,...,2^{k},\quad k=0,...,m.
\end{equation*}
Since $\Delta _{k,j}=\Delta _{k+1,2j-1}+\Delta _{k+1,2j},$ it is clear that $%
I_{k,j}^{m}=I_{k+1,2j-1}^{m}+I_{k+1,2j-1}^{m},$ for $j=1,...,2^{k}$. In
particular $J_{M}=I_{0,1}^{M}=\{1,...,n\}.$ We leave to the reader to show
that each set $I_{k,j}^{m}$ contains at least two elements, if the constant $%
n_{\min }$ is large enough.

\begin{proposition}
\label{Prop-IMKJ}Assume that $n_{\min }>2C_{\max }/C_{\min }\geq 2.$ Then
for any $j=1,...,2^{k},$ $k=0,...,m,$ $m=0,...,M,$ we have $%
\#I_{k,j}^{m}\geq 2.$
\end{proposition}

In the sequel we shall assume that $n>n_{\min }\geq 2C_{\max }/C_{\min }\geq
2.$ Now the sequence $X_{i}^{m},$ $i\in J_{m}$ can be split into dyadic
blocks corresponding to the sets of indices $I_{k,j}^{m}$ as follows: 
\begin{equation*}
\{X_{i}^{m}:i\in J_{m}\}=\sum_{j=1}^{2^{k}}\{X_{i}:i\in I_{k,j}^{m}\},\quad
k=0,...,m.
\end{equation*}
Set 
\begin{equation}
X_{k,j}^{m}=\sum_{i\in I_{k.j}^{m}}X_{i}^{m},\quad B_{k,j}^{m}=E^{\prime
}(X_{k,j}^{m})^{2}=\sum_{i\in I_{k.j}^{m}}E^{\prime }(X_{i}^{m})^{2}.
\label{BCN.6}
\end{equation}
The following assertions are crucial in the proof of our results, as shall
be seen later. The proofs being elementary are left to the reader.

\begin{proposition}
\label{P B1}For any $k=0,...,M-1$ and $j=1,...,2^{k}$ we have 
\begin{equation}
\left| B_{k+1,2j-1}^{m}-B_{k+1,2j}^{m}\right| \leq c\lambda _{n}^{2}.
\label{BCN-DIFF-1}
\end{equation}
\end{proposition}

\begin{proposition}
\label{P B2}For any $k=0,...,M-1$ and $j=1,...,2^{k}$ we have 
\begin{equation*}
c^{-1}\leq B_{k+1,2j-1}^{m}/B_{k+1,2j}^{m}\leq c.
\end{equation*}
\end{proposition}

\subsection{The construction\label{sec: CONSTR}}

Recall that at this moment we are given just two sequences of independent
r.v.'s: $X_{i},$ $i=1,...,n$ on the probability space $(\Omega ^{\prime },%
\mathcal{F}^{\prime },P^{\prime })$ and $N_{i},$ $i=1,...,n$ on the
probability space $(\Omega ,\mathcal{F},P).$ We would like to construct a
sequence of independent r.v.'s $\widetilde{X}_{i},$ $i=1,...,n$ on the
probability space $(\Omega ,\mathcal{F},P)$ such that each $\widetilde{X}%
_{i} $ has the same distribution as $X_{i}$ and the two sequences $%
\widetilde{X}_{i},$ $i=1,...,n$ and $N_{i},$ $i=1,...,n$ are as close as
possible. Before proceeding with the construction we shall describe two
necessary ingredients: the dyadic scheme of Koml\'{o}s, Major and
Tusn\'{a}dy \cite{KMT1} and an auxiliary construction.

\subsubsection{The Koml\'{o}s-Major-Tusn\'{a}dy dyadic scheme \label%
{sec:DYAD-PROC}}

In this section we shall describe a version of the construction appropriate
for our purposes.

Let $\xi _{m,j},$ $j=1,...,2^{m}$ be a sequence of r.v.'s of zero means and
finite variances given on a probability space $(\Omega ^{\prime },\mathcal{F}%
^{\prime },P^{\prime })$, and let $\eta _{m,j}$ $j=1,...,2^{m}$ be a
sequence of normal r.v.'s with the same means and variances given on a
probability space $(\Omega ,\mathcal{F},P).$ At this moment it is not
necessary to assume that these are sequences of independent r.v.'s. The goal
is to construct a version of $\xi _{m,j},$ $j=1,...,2^{m}$ on the
probability space $(\Omega ,\mathcal{F},P).$ The new sequence will be
denoted $\widetilde{\xi }_{m,j},$ $j=1,...,2^{m}$.

Set $\xi _{k,j}=\xi _{k+1,2j-1}+\xi _{k+1,2j}$ and $\eta _{k,j}=\eta
_{k+1,2j-1}+\eta _{k+1,2j},$ for $j=1,...,2^{k}$ and $k=0,...,m-1.$ First
define $\widetilde{\xi }_{0,1}$ to be the quantile transformation of $\eta
_{0,1},$ i.e. define $\widetilde{\xi }_{0,1}$ to be the solution of the
equation 
\begin{equation*}
F_{\xi _{0,1}}\left( \widetilde{\xi }_{0,1}\right) =\Phi _{\eta
_{0,1}^{{}}}\left( \eta _{0,1}\right)
\end{equation*}
where $F_{\xi _{0,1}}\left( x\right) $ is the distribution function of $\xi
_{0,1},$ and $\Phi _{\eta _{0,1}^{{}}}\left( x\right) $ is the distribution
function of $\eta _{0,1}$ (see Section \ref{sec:QT}). Suppose that for some $%
k=0,...,m-1$ the r.v.'s $\widetilde{\xi }_{k,j},$ $j=1,...,2^{k}$ have
already been constructed, and the goal is to construct $\widetilde{\xi }%
_{k+1,j},$ $j=1,...,2^{k+1}$. To this end set for $j=1,...,2^{k}$%
\begin{equation}
V_{k,j}=\alpha _{k+1,2j}\eta _{k+1,2j-1}-\alpha _{k+1,2j-1}\eta _{k+1,2j},
\label{VKJ-not}
\end{equation}
where 
\begin{equation*}
\alpha _{k+1,2j-1}=\left( \frac{B_{k+1,2j-1}}{B_{k+1,2j}}\right)
^{1/2},\quad \alpha _{k+1,2j}=\left( \frac{B_{k+1,2j}}{B_{k+1,2j-1}}\right)
^{1/2}
\end{equation*}
and 
\begin{equation*}
B_{k+1,2j-1}=E\xi _{k+1,2j-1}^{2},\quad B_{k+1,2j}=E\xi _{k+1,2j}^{2}.
\end{equation*}
Define $\widetilde{T}_{k,j}$ to be the conditional quantile transformation
of $V_{k,j}$ given $\widetilde{\xi }_{k,j}$, i.e. for $j=1,...,2^{k}$ define 
$\widetilde{T}_{k,j}$ as the solution of the equation 
\begin{equation}
F_{T_{k,j}|\xi _{k,j}}\left( \widetilde{T}_{k,j}|\widetilde{\xi }%
_{k,j}\right) =\Phi _{V_{k,j}}\left( V_{k,j}\right)  \label{XX-0}
\end{equation}
where $F_{T_{k,j}|\xi _{k,j}}\left( x|y\right) $ is the conditional
distribution function of $T_{k,j}$ given $\xi _{k,j}=y$, and $\Phi
_{V_{k,j}}(x)$ is the distribution function of $V_{k,j}$ (see Section \ref
{sec:QT}). For any $j=1,...,2^{k}$, the desired r.v.'s $\widetilde{\xi }%
_{k+1,2j-1}$ and $\widetilde{\xi }_{k+1,2j}$ are defined as the solution the
linear system 
\begin{equation}
\left\{ 
\begin{array}{l}
\widetilde{T}_{k,j}=\alpha _{k+1,2j}\widetilde{\xi }_{k+1,2j-1}-\alpha
_{k+1,2j-1}\widetilde{\xi }_{k+1,2j}, \\ 
\widetilde{\xi }_{k,j}=\widetilde{\xi }_{k+1,2j-1}+\widetilde{\xi }_{k+1,2j},
\end{array}
\right.  \label{SYST-dyad}
\end{equation}
the determinant of which is obviously strictly positive. This completes
description of the dyadic procedure.

The following result concerns basic properties of the resulting sequence $%
\widetilde{\xi }_{m,j},$ $j=1,...,2^{m}.$

\begin{lemma}
\label{L B1}Assume that $\xi _{m,j},$ $j=1,...,2^{m},$ and $\eta _{m,j},$ $%
j=1,...,2^{m}$ are sequences of independent r.v.'s. Then for any $k=0,...,m,$
the r.v.'s $\widetilde{\xi }_{k,j},$ $j=1,...,2^{k}$ are independent and
such that $\widetilde{\xi }_{k,j}\overset{d}{=}\xi _{k,j},$ $j=1,...,2^{k}$.
Moreover $\widetilde{\xi }_{k,j},$ $j=1,...,2^{k}$ are functions of the
sequence $\eta _{k,j},$ $j=1,...,2^{k}$ only.
\end{lemma}

\begin{proof}
The proof is similar to statements in Koml\'{o}s, Major and Tusn\'{a}dy \cite
{KMT1} (see also Sakhanenko \cite{Sakh}, Einmahl \cite{Einm}, Zaitsev \cite
{Zaits2001}) and therefore will not be detailed here.
\end{proof}

It turns out that the properties of the Koml\'{o}s-Major-Tusn\'{a}dy dyadic
construction established in Lemma \ref{L B1} are sufficient for proving a
strong approximation result if the index functions of the process belong to
the class of indicators. However for proving our functional version we need
one more property of this construction, which we formulate below. Recall
that $V_{k,j}$ and $\widetilde{T}_{k,j}$ are defined by(\ref{VKJ-not}) and (%
\ref{SYST-dyad}).

\begin{lemma}
\label{L B2}If $\xi _{m,j},$ $j=1,...,2^{m},$ and $\eta _{m,j},$ $%
j=1,...,2^{m},$ are sequences of independent r.v.'s, then, for any $%
k=0,...,m,$ the r.v.'s $\widetilde{T}_{k,j}-V_{k,j},$ $j=1,...,2^{k},$ are
independent.
\end{lemma}

\begin{proof}
For the proof of this statement it suffces to note that for any $k=0,...,m,$%
\begin{equation*}
\left\{ \widetilde{\xi }_{k,j},V_{k,j}:j=1,...,2^{k}\right\} 
\end{equation*}
is a collection of of jointly independent random variables.
\end{proof}

\subsubsection{An auxiliary construction\label{sec:AUX-CON}}

In the sequel we shall need also an auxiliary procedure which is not as
powerful as the KMT construction, but which permits us to construct somehow
the components inside an already constructed arbitrary sum of independent
r.v.'s. Below we present one of the possible methods.

We start from an arbitrary sequence of r.v.'s $\xi _{1},...,\xi _{n}$ (not
necessarily independent) given on $(\Omega ^{\prime },\mathcal{F}^{\prime
},P^{\prime })$. Set $S_{k}=\xi _{1}+...+\xi _{k},$ $k=1,...,n.$ Suppose
that on another probability space $(\Omega ,\mathcal{F},P)$ we have
constructed only the r.v.\ $\widetilde{S}_{n}\overset{d}{=}S_{n},$ which
corresponds to the sum $S_{n}$ and we wish to construct its components, i.e. 
$\widetilde{\xi }_{1},...,\widetilde{\xi }_{n}$ such that $\widetilde{\xi }%
_{1}\overset{d}{=}\xi _{1},...,\widetilde{\xi }_{n}\overset{d}{=}\xi _{n}$
and $\widetilde{S}_{n}=\widetilde{\xi }_{1}+...+\widetilde{\xi }_{n}.$ As a
prerequisite we assume that on the probability space $(\Omega ,\mathcal{F}%
,P) $ we are given a sequence of nondegenerate normal r.v.'s $\eta
_{1},...,\eta _{n}$ (not necessarily independent). First we define $%
\widetilde{\xi }_{n}$ to be the conditional quantile transformation of $\eta
_{n}$ given $\widetilde{S}_{n}$, i.e. we define $\widetilde{\xi }_{n}$ to be
the solution of the equation 
\begin{equation*}
F_{\xi _{n}|S_{n}}\left( \widetilde{\xi }_{n}|\widetilde{S}_{n}\right) =\Phi
_{\eta _{n}}\left( \eta _{n}\right)
\end{equation*}
where $F_{\xi _{n}|S_{n}}(x|y)$ is the conditional distribution of $\xi _{n}$
given $S_{n}=y,$ and $\Phi _{\eta _{n}}(x)$ is the distribution function of $%
\eta _{n}.$ Set $\widetilde{S}_{n-1}=\widetilde{S}_{n}-\widetilde{\xi }_{n}.$
If for some $2\leq k\leq n-1$ the r.v.'s\ $\widetilde{\xi }_{n},...,%
\widetilde{\xi }_{k+1}$ and $\widetilde{S}_{k}$ are already constructed, we
define $\widetilde{\xi }_{k}$ to be the conditional quantile transformation
of $\eta _{k}$ given $\widetilde{S}_{k},$ i.e. we define $\ \widetilde{\xi }%
_{k}$ to be the solution of the equation 
\begin{equation*}
F_{\xi _{k}|S_{k}}\left( \widetilde{\xi }_{k}|\widetilde{S}_{k}\right) =\Phi
_{\eta _{k}}\left( \eta _{k}\right)
\end{equation*}
where $F_{\xi _{k}|S_{k}}(x|y)$ is the conditional distribution of $\xi _{k}$
given $S_{k}=y,$ and $\Phi _{\eta _{k}}(x)$ is the distribution function of $%
\eta _{k}$. Set $\widetilde{S}_{k-1}=\widetilde{S}_{k}-\widetilde{\xi }_{k}$%
. Finally, for $k=1,$ we define $\widetilde{\xi }_{1}=\widetilde{S}_{1},$
this completing our procedure.

The easy proof of the following assertion is left to the reader.

\begin{lemma}
\label{Lemma AUX-CNSTR}Assume that $\xi _{1},...,\xi _{n}$ and $\eta
_{1},...,\eta _{n}$ are sequences of independent r.v.'s. Then $\widetilde{%
\xi }_{1},...,\widetilde{\xi }_{n}$ are independent, $\widetilde{\xi }_{i}%
\overset{d}{=}\xi _{i},$ $i=1,...,n,$ and $\widetilde{\xi }_{1}+...+%
\widetilde{\xi }_{n}=\widetilde{S}_{n}.$ Moreover $\widetilde{\xi }_{1},...,%
\widetilde{\xi }_{n}$ are functions of $\eta _{1},...,\eta _{n}$ and $%
\widetilde{S}_{n}$ only.
\end{lemma}

\subsubsection{The main construction}

Our next step is to describe a construction which will result in the desired
sequence $\widetilde{X}_{i},$ $i=1,...,n$. It should be noted that although
both the dyadic procedure and the auxiliary construction described above
work with arbitrary distributions, in order to use the quantile inequalities
stated in Section \ref{sec:QT} (which actually will provide the desired
closeness of $\widetilde{X}_{i},$ $i=1,...,n,$ and $N_{i},$ $i=1,...,n$),
one has to assume the r.v.'s $X_{i},$ $i=1,...,n$ to be in the class $%
\mathfrak{D}_{0}(r),$ for some $r>0,$ or to be identically distributed (as
in Koml\'{o}s, Major and Tusn\'{a}dy \cite{KMT1}, \cite{KMT2}). In order to
avoid such assumptions we shall employ an inductive procedure which goes
back to the paper of Sakhanenko \cite{Sakh}. The idea is first to substitute
the initial sequence with some smoothed sequences, and then to apply the
dyadic procedure described in Section \ref{sec:DYAD-PROC} to the smoothed
sequences. Below we formally describe this construction.

Consider the product probability space $(\Omega ^{\prime \prime },\mathcal{F}%
^{\prime \prime },P^{\prime \prime })=(\Omega ^{\prime },\mathcal{F}^{\prime
},P^{\prime })\times (\Omega ,\mathcal{F},P),$ where $P^{\prime \prime
}=P^{\prime }\times P.$ It is obvious that the sequences $X_{i},$ $%
i=1,...,n, $ and $N_{i},$ $i=1,...,n$ are independent on the probability
space $(\Omega ^{\prime \prime },\mathcal{F}^{\prime \prime },P^{\prime
\prime })$.

Recall that above we introduced the sets of indices $J_{m}=\{i:1\leq
i2^{M-m}\leq n\}.$ For each $m=M,...,0,$ the set $J_{m}$ can be decomposed
as $J_{m}=J_{m}^{1}+J_{m}^{2},$ where 
\begin{equation*}
J_{m}^{1}=\{i\text{-odd}:i\in J_{m}\},\quad J_{m}^{2}=\{i\text{-even}:i\in
J_{m}\},\quad m=M,...,1,
\end{equation*}
and $J_{0}^{1}=J_{0},$ $J_{0}^{2}=\emptyset .$ It is clear that 
\begin{equation*}
J_{m-1}=\{i:2i\in J_{m}\},\quad m=1,...,M.
\end{equation*}
To start our iterative construction, for any $i\in J_{M}=\{1,...,n\},$
define the following r.v.'s: 
\begin{equation}
X_{2i}^{M+1}=X_{i},\quad \widetilde{Y}_{2i}^{M+1}=N_{i}.  \label{BC.1}
\end{equation}
We proceed to describe the $m$-th step of our construction which is
performed consecutively for all $m=M,...,0.$

$\bullet \,\,m$\textit{-th step.} For any $i\in J_{m},$ define the following
r.v.'s: 
\begin{equation}
X_{i}^{m}=X_{2i}^{m+1},\quad W_{i}^{m}=\widetilde{Y}_{2i}^{m+1},
\label{BC.4}
\end{equation}
and 
\begin{equation}
Y_{i}^{m}=\left\{ 
\begin{array}{ll}
X_{i}^{m}, & \quad \text{if}\quad i\in J_{m}^{1}, \\ 
W_{i}^{m}, & \quad \text{if}\quad i\in J_{m}^{2}.
\end{array}
\right.  \label{BC.4a}
\end{equation}
Note that the r.v.'s $W_{i}^{m},$ $i\in J_{m},$ are defined on the
probability space $(\Omega ,\mathcal{F},P),$ while the r.v.'s $Y_{i}^{m},$ $%
i\in J_{m},$ are defined on the probability space $(\Omega ^{\prime \prime },%
\mathcal{F}^{\prime \prime },P^{\prime \prime }).$ Here $X_{i}^{m},$ $i\in
J_{m}$ is the part of the initial sequence $X_{i},$ $i\in J_{M}=\{1,...,n\}$
given on the probability space $(\Omega ^{\prime },\mathcal{F}^{\prime
},P^{\prime }),$ which is not yet constructed on the probability space $%
(\Omega ,\mathcal{F},P)$; $W_{i}^{m},$ $i\in J_{m}$ is the corresponding
sequence of normal r.v.'s and $Y_{i}^{m},$ $i\in J_{m}$ is the smoothed
sequence which is constructed at this step. Consider the following sums: for 
$j=1,...,2^{k}$ and $k=0,...,m$ set 
\begin{equation}
Y_{k,j}^{m}=\sum_{i\in I_{k,j}^{m}}Y_{i}^{m},\quad W_{k,j}^{m}=\sum_{i\in
I_{k,j}^{m}}W_{i}^{m}.  \label{Ymmkj-def}
\end{equation}
Then obviously for $j=1,...,2^{k}$ and $k=0,...,m-1$%
\begin{equation}
Y_{k,j}^{m}=Y_{k+1,2j-1}^{m}+Y_{k+1,2j}^{m},\quad
W_{k,j}^{m}=W_{k+1,2j-1}^{m}+W_{k+1,2j}^{m}.  \notag
\end{equation}
We will apply the dyadic procedure described in Section \ref{sec:DYAD-PROC},
with $\xi _{m,j}=Y_{m,j}^{m}$ and $\eta _{m,j}=W_{m,j}^{m},$ $j=1,...,2^{m},$
to construct a doubly indexed sequence $\widetilde{Y}_{kj}^{m},$ $%
j=1,...,2^{m},$ $k=0,...,m.$ Let $\widetilde{Y}_{0,1}^{m}$ be the quantile
transformation of $W_{0,1}^{m},$ i.e. let $\widetilde{Y}_{0,1}^{m}$ be the
solution of the equation 
\begin{equation}
F_{Y_{0,1}^{m}}\left( \widetilde{Y}_{0,1}^{m}\right) =\Phi
_{W_{0,1}^{m}}\left( W_{0,1}^{m}\right) ,  \label{CNS-mm01}
\end{equation}
where $F_{Y_{0,1}^{m}}\left( x\right) $ is the distribution function of $%
Y_{0,1}^{m}$ and $\Phi _{W_{0,1}^{m}}\left( x\right) $ is the distribution
function of $W_{0,1}^{m}.$ The solution exists since $W_{0,1}^{m}$ is a
nondegenerate normal r.v. Assume that we have already constructed $%
\widetilde{Y}_{k,j}^{m},$ $j=1,...,2^{k},$ for some $k=0,...,m-1.$ We shall
construct such an array with $k+1$ replacing $k.$ To this end set, for $%
j=1,...,2^{k},$%
\begin{equation}
V_{k,j}^{m}=\alpha _{k+1,2j}^{m}W_{k+1,2j-1}^{m}-\alpha
_{k+1,2j-1}^{m}W_{k+1,2j}^{m},  \label{BC-2-m}
\end{equation}
where 
\begin{equation*}
\alpha _{k+1,2j-1}^{m}=\left( \frac{B_{k+1,2j-1}^{m}}{B_{k+1,2j}^{m}}\right)
^{1/2},\quad \alpha _{k+1,2j}^{m}=\left( \frac{B_{k+1,2j}^{m}}{%
B_{k+1,2j-1}^{m}}\right) ^{1/2}
\end{equation*}
and 
\begin{equation}
B_{k+1,2j-1}^{m}=E\left( \widetilde{Y}_{k+1,2j-1}^{m}\right) ^{2},\quad
B_{k+1,2j}^{m}=E\left( \widetilde{Y}_{k+1,2j}^{m}\right) ^{2}.  \label{BBB-m}
\end{equation}
Let $\widetilde{T}_{k,j}^{m}$ be the conditional quantile transformation of $%
V_{k,j}^{m},$ given $\widetilde{Y}_{k,j}^{m},$ for $j=1,...,2^{k},$ i. e.
let $\widetilde{T}_{k,j}^{m}$ be the solution of the equation 
\begin{equation}
F_{T_{k,j}^{m}|Y_{k,j}^{m}}\left( \widetilde{T}_{k,j}^{m}|\widetilde{Y}%
_{k,j}^{m}\right) =\Phi _{V_{k,j}^{m}}\left( V_{k,j}^{m}\right) ,
\label{CNS-mmkj}
\end{equation}
where $F_{T_{k,j}^{m}|Y_{k,j}^{m}}\left( x|y\right) $ is the conditional
distribution function of $T_{k,j}^{m},$ given $Y_{k,j}^{m},$ and $\Phi
_{W_{k,j}^{m}}\left( x\right) $ is the distribution function of $%
W_{k,j}^{m}. $ The solution exists, since $V_{k,j}^{M}$ is a nondegenerate
normal r.v. For any $j=1,...,2^{k}$ we define the desired r.v.'s $\widetilde{%
Y}_{k+1,2j-1}^{m}$ and $\widetilde{Y}_{k+1,2j}^{m}$ as the solution of the
linear system 
\begin{equation}
\left\{ 
\begin{array}{l}
\widetilde{T}_{k,j}^{m}=\alpha _{k+1,2j}^{m}\widetilde{Y}_{k+1,2j-1}^{m}-%
\alpha _{k+1,2j-1}^{m}\widetilde{Y}_{k+1,2j}^{m}, \\ 
\widetilde{Y}_{k,j}^{m}=\widetilde{Y}_{k+1,2j-1}^{m}+\widetilde{Y}%
_{k+1,2j}^{m}.
\end{array}
\right.  \label{SYST-mm}
\end{equation}
Thus the r.v.'s $\widetilde{Y}_{k,j}^{m},$ $j=1,...,2^{k}$ are constructed
for all $k=0,...,m$ on the probability space $(\Omega ,\mathcal{F},P)$. It
remains to construct the components inside each sum $\widetilde{Y}%
_{m,j}^{m}, $ $j=1,...,2^{m}.$ For this we make use of the auxiliary
construction described in Section \ref{sec:AUX-CON}, with $\xi _{i}\equiv
Y_{i}^{m}$ and $\eta _{i}\equiv W_{i}^{m},$ $i\in I_{m,j}^{m}.$ For each
fixed $j$ and $m$ it provides a sequence of r.v.'s $\widetilde{Y}%
_{i}^{m}\equiv \widetilde{\xi }_{i},$ $i\in I_{m,j}^{m},$ such that 
\begin{equation}
\widetilde{Y}_{m,j}^{m}=\sum_{i\in I_{m,j}^{m}}\widetilde{Y}_{i}^{m}.
\label{YWmm-sum}
\end{equation}
This completes the $m$-th step of our construction.

Let us recall briefly some notation associated with the construction, which
will also be used in the sequel. For any $m=M,...,0$ we have defined the
r.v.'s $Y_{i}^{m},$ $W_{i}^{m},$ $\widetilde{Y}_{i}^{m},$ $i\in J_{m},$ and $%
Y_{k,j}^{m},$ $W_{k,j}^{m},$ $\widetilde{Y}_{k,j}^{m},$ $j=1,...,2^{k},$ $%
k=0,...,m,$ such that, by (\ref{Ymmkj-def}) and (\ref{YWmm-sum}) (cp. with (%
\ref{BCN.6})), 
\begin{equation}
Y_{k,j}^{m}=\sum_{i\in I_{k,j}^{m}}Y_{i}^{m},\quad W_{k,j}^{m}=\sum_{i\in
I_{k,j}^{m}}W_{i}^{m},\quad \widetilde{Y}_{k,j}^{m}=\sum_{i\in I_{k,j}^{m}}%
\widetilde{Y}_{i}^{m},  \label{YYY}
\end{equation}
for $k=0,...,m,$ $j=1,...,2^{k},$ $m=0,...,M.$

\subsection{Correctness and some useful properties}

In fact implicitly the construction of the desired sequence $\widetilde{X}%
_{i},$ $i=1,...,n$ has already been carried out; it remains to select the
appropriate components from the sequences $\{\widetilde{Y}_{i}^{m}:i\in
J_{m}\}$ found above. But before this step we need to show that the
construction is performed correctly, and we shall also discuss some
properties of the r.v.'s $\widetilde{Y}_{i}^{m}$ and $W_{i}^{m}$ introduced.
The proofs of the following assertions are left to the reader.

In analogy to $X_{i}^{m}$ (see (\ref{Xmi})), set $N_{i}^{m}=N_{i2^{M-m}},$
where $m=0,...,M,$ $i\in J_{m}.$

\begin{lemma}
\label{Lemma-CNSTR}For any $m=0,...,M$ the following statements hold true:

a) The r.v.'s $W_{i}^{m},$ $i\in J_{m},$ are independent and satisfy $%
W_{i}^{m}\overset{d}{=}N_{i}^{m}=N_{i2^{M-m}},$ $i\in J_{m}.$

b) The r.v.'s $\widetilde{Y}_{i}^{m},$ $i\in J_{m}$ are independent, are
functions of $W_{i}^{m},$ $i\in J_{m}$ only and satisfy, for $i\in J_{m},$ 
\begin{equation*}
\widetilde{Y}_{i}^{m}\overset{d}{=}Y_{i}^{m}\overset{d}{=}\left\{ 
\begin{array}{ll}
X_{i}^{m} & \quad \text{if}\quad i\in J_{m}^{1}, \\ 
N_{i}^{m} & \quad \text{if}\quad i\in J_{m}^{2}.
\end{array}
\right. 
\end{equation*}
\end{lemma}

\begin{remark}
Since by Proposition \ref{Prop-IMKJ} $\#I_{k,j}^{m}\geq 2,$ from Lemma \ref
{Lemma-CNSTR} and from (\ref{YYY}) it follows that $W_{k,j}^{m},$ $%
j=1,...,2^{k},$ $k=0,...,m$ are nondegenerate normal r.v.'s which ensures
that the solutions of the equations (\ref{CNS-mm01}), (\ref{CNS-mmkj})
exist. This proves the correctness of the main construction.
\end{remark}

\begin{proposition}
\label{P Q3}The vectors $\left\{ \widetilde{Y}_{i}^{m}:i\in
J_{m}^{1}\right\} ,$ $m=M,...,0$ are independent.
\end{proposition}

Now finally we are able to present the sequence $\widetilde{X}_{i},$ $%
i=1,...,n.$ It is defined on the probability space $(\Omega ,\mathcal{F},P)$
in the following way: 
\begin{equation}
\widetilde{X}_{i2^{M-m}}=\widetilde{Y}_{i}^{m},\quad \text{where}\quad i\in
J_{m}^{1},\quad 0\leq m\leq M.  \label{BC.d1}
\end{equation}

\begin{proposition}
\label{P BC1} $\widetilde{X}_{i},$ $i=1,...,n,$ are independent and such
that $\widetilde{X}_{i}\overset{d}{=}X_{i},$ $i=1,...,n.$
\end{proposition}

\begin{proof}
The required assertion follows from Lemma \ref{Lemma-CNSTR} and Proposition 
\ref{P Q3}.
\end{proof}

In the proof of our main result Theorem \ref{T3-2}, the following elementary
representation is essential. Recall that $t_{i}^{m}=t_{\nu }=\nu /n$ where $%
\nu =i2^{M-m},$ $i\in J_{m},$ $m=0,...,M$ (see Section \ref{sec:SNot}).

\begin{proposition}
\label{P R1}For any real valued function $f(t)$ on the interval $[0,1],$ we
have 
\begin{equation*}
\sum_{i=1}^{n}f(t_{i})\left( \widetilde{X}_{i}-N_{i}\right)
=\sum_{m=0}^{M}\sum_{j\in J_{m}}f(t_{j}^{m})\left( \widetilde{Y}%
_{j}^{m}-W_{j}^{m}\right) .
\end{equation*}
\end{proposition}

\subsection{Quantile inequalities\label{sec:QI}}

In this section we shall establish so-called quantile inequalities (see
Lemma \ref{L B4} and Lemma \ref{L B5}), which will ensure the required
closeness of the r.v.'s $\widetilde{X}_{i},$ $i=1,...,n,$ and $N_{i},$ $%
i=1,...,n.$

The following lemma shows that the r.v.'s $Y_{k,j}^{m},$ $j=1,...,2^{k}$ are
smooth enough to allow application of the quantile inequalities stated in
Section \ref{sec:QT}.

\begin{lemma}
\label{L B3}For $m=0,...,M,$ $k=0,...,m,$ $j=1,...,2^{k}$ the r.v.\ $%
Y_{k,j}^{m}$ is in the class $\mathfrak{D}(r),$ for some positive absolute
constant $r.$
\end{lemma}

\begin{proof}
We shall check conditions (\ref{Q.1}), (\ref{Q.2}) and (\ref{Q.3}) in
Section \ref{sec:QT}. Toward this end fix $m,$ $k,$ $j$ as in the condition
of the lemma and note that 
\begin{equation*}
\zeta _{0}\equiv Y_{k,j}^{m}=\sum_{i\in I_{k,j}^{m}}Y_{i}^{m}=\sum_{i\in
I_{1}}Y_{i}^{m}+\sum_{i\in I_{2}}Y_{i}^{m}\equiv \zeta _{1}+\zeta _{2},
\end{equation*}
where $I_{1}$ and $I_{2}$ are the sets of all odd and even indices in $%
I_{k,j}^{m}$ respectively. By Lemma \ref{Lemma-CNSTR}, we have $Y_{i}^{m}%
\overset{d}{=}N_{i},$ for any $i\in I_{2}.$ Thus \ $\zeta _{2}$ is actually
a sum of independent normal r.v.'s. Since $n_{\min }$ is large enough, the
set $I_{k,j}^{m}$ has at least two elements (see Proposition \ref{Prop-IMKJ}%
), from which we conclude that $I_{2}$ has at least one element. Next,
taking into account (\ref{NR.1}) and the obvious inequality $\#I_{2}\geq 
\frac{1}{3}\#I_{k,j}^{m},$ we get 
\begin{equation*}
E\zeta _{2}^{2}\geq C_{\min }\lambda _{n}\#I_{2}\geq \frac{C_{\min }}{3}%
\lambda _{n}\#I_{k,j}^{m}\geq cE\zeta _{0}^{2}.
\end{equation*}
For $|h|\leq \lambda $ and $t\in R,$ let 
\begin{equation*}
f_{\zeta _{i}}(t,h)=E\exp \{(\mathbf{i}t+h)\zeta _{i}\}/E\exp \{h\zeta _{i}\}
\end{equation*}
be the conjugate characteristic function of the r.v.\ $\zeta _{i},$ $i=0,1,2.
$ Since $\zeta _{1}$ and $\zeta _{2}$ are independent and $\zeta _{2}$ is
normal, 
\begin{eqnarray*}
|f_{\zeta _{0}}(t,h)| &=&|f_{\zeta _{1}}(t,h)\,\,f_{\zeta _{2}}(t,h)|\leq
|f_{\zeta _{2}}(t,h)| \\
\  &\leq &\exp \{-\frac{t^{2}}{2}E\zeta _{2}^{2}\}\leq \exp \{-\frac{t^{2}}{2%
}cE\zeta _{0}^{2}\},
\end{eqnarray*}
for $|h|\leq \lambda ,$ $t\in \mathbf{R}^{1}.$ With this bound we have 
\begin{equation*}
\int_{|t|>\varepsilon }|f_{\zeta _{0},h}(t)|dt\leq \int_{|t|>\varepsilon
}\exp \{-\frac{t^{2}}{2}cE\zeta _{0}^{2}\}dt\leq \frac{\mu }{\varepsilon
E\zeta _{0}^{2}},
\end{equation*}
where $\mu $ is some absolute constant, which proves that $\zeta
_{0}=Y_{k,j}^{m}$ satisfies condition (\ref{Q.3}).

It remains only to show that conditions (\ref{Q.1}) and (\ref{Q.2}) are
satisfied. The first condition follows from (\ref{NR.1}) as soon as $%
Y_{i}^{m}\overset{d}{=}X_{i}$ or $Y_{i}^{m}\overset{d}{=}N_{i}$ for any $%
i\in I_{k,j}^{m}\subseteq J_{m},$ by Lemma \ref{Lemma-CNSTR}. For the second
we make use of (\ref{NR.2}) and of the elementary fact that Sakhanenko's
condition (\ref{Q.2}) holds true for any normal r.v.\ $N$ if $\lambda $ is
small enough: $\lambda \leq c(\mathrm{Var}N)^{-1/2}$ (see Remark \ref
{Rem-S-for-norm}).
\end{proof}

Recall that for any $m=0,...,M$, $k=1,...,m$ and $j=1,...,2^{k},$ by (\ref
{SYST-mm}), 
\begin{equation}
\widetilde{T}_{k,j}^{m}=\alpha _{k,2j}^{m}\widetilde{Y}_{k,2j-1}^{m}-\alpha
_{k,2j-1}^{m}\widetilde{Y}_{k,2j}^{m},  \label{QI.1x-1}
\end{equation}
and by (\ref{BC-2-m}), 
\begin{equation}
V_{k,j}^{m}=\alpha _{k,2j}^{m}W_{k,2j-1}^{m}-\alpha
_{k,2j-1}^{m}W_{k,2j}^{m}.  \label{QI.1x-2}
\end{equation}
Recall also that $B_{k,j}^{m}=E^{\prime }(X_{k,j}^{m})^{2}$ (see (\ref{BCN.6}%
)).

The following quantile inequalities show that the r.v.'s $\widetilde{T}%
_{k,j}^{m}$ and $W_{k,j}^{m}$ are close enough. These statements are crucial
for our results.

\begin{lemma}
\label{L B4}For any $m=0,...,M,$ we have 
\begin{equation*}
\left| \widetilde{Y}_{0,1}^{m}-W_{0,1}^{m}\right| \leq c_{1}\left\{ 1+\frac{%
\left( \widetilde{Y}_{0,1}^{m}\right) ^{2}}{B_{0,1}^{m}}\right\} ,
\end{equation*}
provided $\left| \widetilde{Y}_{0,1}^{m}\right| \leq c_{2}B_{0,1}^{m}$ and $%
B_{0,1}^{m}\geq c_{3},$ where $c_{1},$ $c_{2}$ and $c_{3}$ are positive
absolute constants.
\end{lemma}

\begin{proof}
According to the construction, $\widetilde{Y}_{0,1}^{m}$ is the quantile
transformation of $W_{0,1}^{m}$ (see (\ref{CNS-mm01})). Then it suffices to
note that, by Lemma \ref{L B3}, the r.v.\ $\widetilde{Y}_{0,1}^{m}$ is in
the class $\mathfrak{D}(\lambda _{0})$ and to apply Lemma \ref{L Q1} with $%
X=Y_{0,1}^{m},$ $N=W_{0,1}^{m}$ and $\widetilde{X}=\widetilde{Y}_{0,1}^{m}.$
\end{proof}

\begin{lemma}
\label{L B5}Let $m=0,...,M,$ $k=0,...,m-1,$ $j=1,...,2^{k}.$ Then 
\begin{equation*}
\left| \widetilde{T}_{k,j}^{m}-V_{k,j}^{m}\right| \leq c_{1}\left\{ 1+\frac{%
\left( \widetilde{Y}_{k+1,2j-1}^{m}\right) ^{2}}{B_{k+1,2j-1}^{m}}+\frac{%
\left( \widetilde{Y}_{k+1,2j}^{m}\right) ^{2}}{B_{k+1,2j}^{m}}\right\} ,
\end{equation*}
provided $\left| \widetilde{Y}_{k+1,2j-1}^{m}\right| \leq
c_{2}B_{k+1,2j-1}^{m},$ $\left| \widetilde{Y}_{k+1,2j}^{m}\right| \leq
c_{2}B_{k+1,2j}^{m}$ and $B_{k+1,2j-1}^{m}\geq c_{3},$ $B_{k+1,2j}^{m}\geq
c_{3},$ where $c_{1},$ $c_{2}$ and $c_{3}$ are positive absolute constants.
\end{lemma}

\begin{proof}
Fix $m,$ $k,$ and $j$ as in the condition of the lemma. We are going to make
use of Lemma \ref{L Q2} with 
\begin{equation}
\widetilde{X}_{1}=\widetilde{Y}_{k+1,2j-1}^{m},\quad \widetilde{X}_{2}=%
\widetilde{Y}_{k+1,2j}^{m},\quad \widetilde{X}_{0}=\widetilde{X}_{1}+%
\widetilde{X}_{2}=\widetilde{Y}_{k,j}^{m},  \notag
\end{equation}
and 
\begin{equation}
N_{1}=W_{k+1,2j-1}^{m},\quad N_{2}=W_{k+1,2j}^{m},\quad
N_{0}=N_{1}+N_{2}=W_{k,j}^{m}.  \notag
\end{equation}
Note that, by Lemma \ref{L B3}, the r.v.'s $\widetilde{X}_{0},$ $\widetilde{X%
}_{1}$ and $\widetilde{X}_{2}$ are in the class $\mathfrak{D}(r)$ for some
absolute constant $r>0.$ Since by construction $\widetilde{T}_{k,j}^{m}$ is
the conditional quantile transformation of $V_{k,j}^{m}$ (see (\ref{CNS-mmkj}%
)), Lemma \ref{L Q2} implies 
\begin{equation}
\left| \widetilde{T}_{k,j}^{m}-V_{k,j}^{m}\right| \leq c_{1}\frac{B_{0}}{B}%
\left\{ 1+\frac{1}{B^{2}}\left( \widetilde{X}_{1}^{2}+\widetilde{X}%
_{2}^{2}\right) \right\} ,  \label{QI.4}
\end{equation}
provided 
\begin{equation}
\left| \widetilde{T}_{k,j}^{m}\right| \leq c_{2}B^{2},\quad \left| 
\widetilde{Y}_{k,j}^{m}\right| \leq c_{2}B^{2},  \label{QI.10}
\end{equation}
and $B\geq c_{3},$ where 
\begin{equation*}
B_{1}^{2}=B_{k+1,2j-1}^{m},\quad B_{2}^{2}=B_{k+1,2j}^{m},\quad
B_{0}^{2}=B_{1}^{2}+B_{2}^{2},\quad B^{2}=\frac{B_{1}B_{2}}{B_{0}}.
\end{equation*}
By Proposition \ref{P B2}, we have 
\begin{equation}
c_{4}^{-1}\leq B_{1}^{2}/B_{2}^{2}\leq c_{4},.  \label{QI.11-11}
\end{equation}

Now we check that (\ref{QI.10}) holds true if $\left| \widetilde{X}%
_{1}\right| \leq c_{5}B_{1}^{2}$ and $\left| \widetilde{X}_{2}\right| \leq
c_{5}B_{2}^{2},$ where $c_{5}$ is a sufficiently small constant. Indeed 
\begin{equation*}
\left| \widetilde{T}_{k,j}^{m}\right| \leq \frac{B_{2}}{B_{1}}\left| 
\widetilde{X}_{1}\right| +\frac{B_{1}}{B_{2}}\left| \widetilde{X}_{2}\right|
\leq 2c_{5}B_{1}B_{2}.
\end{equation*}
By (\ref{QI.11-11}), we get $\left| \widetilde{T}_{k,j}^{m}\right| \leq
c_{5}c_{6}B^{2}.$ Choosing the constant $c_{6}$ such that $c_{5}c_{6}\leq
c_{2},$ we see that (\ref{QI.10}) is satisfied. Exactly in the same way we
show that the second inequality in (\ref{QI.10}) holds true. Condition $%
B\geq c_{3}$ follows easily from (\ref{QI.11-11}).
\end{proof}

\section{Proof of the main results \label{sec:PMR}}

\subsection{An auxiliary exponential bound}

We keep the same notation as in the previous section. In addition set for
brevity 
\begin{equation}
\widetilde{S}_{0}^{m}=\widetilde{Y}_{0,1}^{m}-W_{0,1}^{m},\quad \widetilde{S}%
_{k,j}^{m}=\widetilde{T}_{k,j}^{m}-V_{k,j}^{m},\quad j=1,...,2^{k},\quad
k=0,...,m,  \label{Aux-st-1}
\end{equation}
where $\widetilde{T}_{k,j}^{m}$ and $V_{k,j}^{m}$ are defined by (\ref
{QI.1x-1}) and (\ref{QI.1x-2}). The main result of this section is Lemma \ref
{lemma-basic} which establishes an exponential type bound for the
differences $\widetilde{S}_{k,j}^{m}$ and $\widetilde{S}_{0}^{m}$. Because
of the special construction of $\widetilde{T}_{k,j}^{m}$ and $V_{k,j}^{m}$
on the same probability space, this bound is much better that the usual
exponential bounds (cf. Lemma \ref{Lemma-exp-bound} below). This statement
plays a crucial role in establishing our functional version of the Hungarian
construction. It is the only place where the quantile inequalities are used.

\begin{lemma}
\label{lemma-basic}For any $m=0,...,M,$ $k=0,...,m-1,$ $j=1,...,2^{k},$%
\begin{equation*}
E\exp \left\{ t\widetilde{S}_{0}^{m}\right\} \leq \exp \left\{
c_{1}t^{2}\right\} ,\quad E\exp \left\{ t\widetilde{S}_{k,j}^{m}\right\}
\leq \exp \left\{ c_{1}t^{2}\right\} ,\quad |t|\leq c_{0}.
\end{equation*}
\end{lemma}

We postpone the proof of the lemma to the end of this section; it will be
based on some estimates stated and proved below.

\begin{lemma}
\label{L P1}For any $\varepsilon >0$ there is a constant $c(\varepsilon )$
depending only on $\varepsilon ,$ such that for any $m=0,...,M,$ $k=0,...,m$
and $j\in J_{k},$%
\begin{equation*}
P\left( \left| \widetilde{Y}_{k,j}^{m}\right| >\varepsilon
B_{k,j}^{m}\right) \leq 2\exp \left\{ -c(\varepsilon )B_{k,j}^{m}\right\} .
\end{equation*}
\end{lemma}

\begin{proof}
By Chebyshev's inequality, we have for $t>0$%
\begin{equation}
P\left( \widetilde{Y}_{k,j}^{m}>\varepsilon B_{k,j}^{m}\right) \leq \exp
\left\{ -t\varepsilon B_{k,j}^{m}\right\} E\exp \{tY_{k,j}^{m}\}.
\label{P.0}
\end{equation}
Note that by (\ref{YYY}) and by Lemma \ref{Lemma-CNSTR}, the r.v. $%
\widetilde{Y}_{k,j}^{m}$ is the sum of independent r.v.'s $\widetilde{Y}%
_{i}^{m},i\in I_{k,j}^{m}.$ Then by (\ref{NR.2}) and Lemma \ref{L A1}, we
obtain for $|t|\leq \lambda /3,$ 
\begin{equation*}
E\exp \left\{ t\widetilde{Y}_{k,j}^{m}\right\} =\prod_{i\in
I_{k,j}^{m}}E\exp \left\{ t\widetilde{Y}_{i}^{m}\right\} \leq \exp \left\{
t^{2}B_{k,j}^{m}\right\} .
\end{equation*}
Inserting this bound into (\ref{P.0}), with an appropriate choice of $t$
(depending on $\varepsilon $), we get 
\begin{equation*}
E\left( \widetilde{Y}_{k,j}^{m}>\varepsilon B_{k,j}^{m}\right) \leq \exp
\left\{ -c(\varepsilon )B_{k,j}^{m}\right\} .
\end{equation*}
In the same way one can show that 
\begin{equation*}
E{(}\widetilde{Y}_{k,j}^{m}<-\varepsilon B_{k,j}^{m})\leq \exp
\{-c(\varepsilon )B_{k,j}^{m}\},
\end{equation*}
which in conjunction with the previous bound proves the lemma.
\end{proof}

\begin{lemma}
\label{Lemma-exp-bound}Let $m=0,...,M,$ $k=0,...,m-1,$ $j=1,...,2^{k}$. Then
for any $0\leq t\leq c_{1}$ we have 
\begin{equation*}
E\exp \left\{ t\left| \widetilde{S}_{k,j}^{m}\right| \right\} \leq c_{2}\exp
\left\{ t^{2}B_{k,j}^{m}\right\} .
\end{equation*}
\end{lemma}

\begin{proof}
Fix $m,$ $k$ and $j$ as in the condition of the lemma. From (\ref{Aux-st-1})
and from the H\"{o}lder inequality one gets, for $0\leq t\leq \lambda ,$ 
\begin{equation}
E\exp \left\{ t|\widetilde{S}_{k,j}^{m}|\right\} \leq \left( E\exp \{t|%
\widetilde{T}_{k,j}^{m}|\}E\exp \{t|V_{k,j}^{m}|\}\right) ^{1/2}.
\label{P.6}
\end{equation}
The r.v.\ $\widetilde{Y}_{k+1,2j-1}^{m}$ and $\widetilde{Y}_{k+1,2j}^{m}$
are independent, hence by (\ref{QI.1x-1}) 
\begin{equation}
E\exp \left\{ t|\widetilde{T}_{k,j}^{m}|\right\} \leq E\exp \left\{ t\alpha
_{k+1,2j}^{m}\left| \widetilde{Y}_{k+1,2j-1}^{m}\right| \right\} E\exp
\left\{ t\alpha _{k+1,2j-1}^{m}\left| \widetilde{Y}_{k+1,2j}^{m}\right|
\right\} .  \label{P.6a}
\end{equation}
Since by (\ref{YYY}) and by Lemma \ref{Lemma-CNSTR}, \ $\widetilde{Y}%
_{k+1,2j-1}^{m}$ is exactly the sum of independent r.v.'s $\widetilde{Y}%
_{i}^{m},$ $i\in I_{k+1,2j-1}^{m},$ one has 
\begin{equation*}
E\exp \left\{ \pm t\alpha _{k+1,2j}^{m}\widetilde{Y}_{k+1,2j-1}^{m}\right\}
=\prod_{i\in I_{k+1,2j-1}^{m}}E\exp \left\{ \pm t\alpha _{k+1,2j}^{m}%
\widetilde{Y}_{i}^{m}\right\} .
\end{equation*}
Taking into account (\ref{NR.2}) and choosing $t$ small enough ($t\leq
\lambda /3$), by Lemma \ref{L A1} one obtains 
\begin{eqnarray*}
E\exp \left\{ \pm t\alpha _{k+1,2j}^{m}\widetilde{Y}_{k+1,2j-1}^{m}\right\} 
&\leq &\prod_{i\in J_{k+1,2j-1}^{m}}E\exp \left\{ t^{2}(\alpha
_{k+1,2j}^{m})^{2}E(X_{i}^{m})^{2}\right\}  \\
&\leq &\exp \left\{ t^{2}(\alpha _{k+1,2j}^{m})^{2}B_{k+1,2j-1}^{m}\right\} .
\end{eqnarray*}
Since $(\alpha _{k+1,2j}^{m})^{2}=B_{k+1,2j}^{m}/B_{k+1,2j-1}^{m},$ 
\begin{equation*}
E\exp \left\{ t\alpha _{k+1,2j}^{m}|\widetilde{Y}_{k+1,2j-1}^{m}|\right\}
\leq 2\exp \left\{ t^{2}B_{k+1,2j}^{m}\right\} .
\end{equation*}
For the second expectation on the right hand side of (\ref{P.6a}) one gets
an analogous bound. Then 
\begin{equation}
E\exp \left\{ t|\widetilde{T}_{k,j}^{m}|\right\} \leq 4\exp \left\{
t^{2}B_{k+1,2j}^{m}+t^{2}B_{k+1,2j-1}^{m}\right\} =4\exp \left\{
t^{2}B_{k,j}^{m}\right\} .  \label{P.7a}
\end{equation}
A similar bound holds for the second expectation on the right-hand side of (%
\ref{P.6}), i. e. 
\begin{equation}
E\exp \left\{ t|V_{k,j}^{m}|\right\} \leq 4\exp \left\{
t^{2}B_{k,j}^{m}\right\} .  \label{P.8}
\end{equation}
Now the lemma follows from (\ref{P.7a}), (\ref{P.8}) and (\ref{P.6}).
\end{proof}

Now we are prepared to show that $\widetilde{S}_{k,j}^m$ has a bounded
exponential moment uniformly in $m,$ $k$ and $j.$

\begin{lemma}
\label{L P2}For any $m=0,...,M,$ $k=0,...,m-1,$ $j=1,...,2^{k}$%
\begin{equation*}
E\exp \left\{ c_{1}\left| \widetilde{S}_{k,j}^{m}\right| \right\} \leq c_{2}.
\end{equation*}
\end{lemma}

\begin{proof}
Fix $m,$ $k$ and $j$ as in the condition of the lemma. It is enough to
consider the case where $B_{k+1,2j-1}^{m}$ and $B_{k+1,2j}^{m}$ are greater
than $c^{\prime }$ only, where $c^{\prime }$ is the absolute constant $c_{3}$
in Lemma \ref{L B5}; otherwise, by Proposition \ref{P B2}, we have $%
B_{k+1,2j-1}^{m},B_{k+1,2j}^{m}\leq c_{1}$ (thus $%
B_{k,j}^{m}=B_{k+1,2j-1}^{m}+B_{k+1,2j}^{m}\leq 2c_{1}$) and the claim
follows from Lemma \ref{Lemma-exp-bound}.

Set for brevity 
\begin{equation}
G_{k+1,l}^{m}=\left\{ \left| \widetilde{Y}_{k+1,l}^{m}\right| \leq c^{\prime
\prime }B_{k+1,l}^{m}\right\} ,\quad l=2j-1,2j,  \label{QI.2x}
\end{equation}
where $c^{\prime \prime }=\min \{1,c_{2}\}$ and $c_{2}$ is the absolute
constant in Lemma \ref{L B5}. Denote by $G_{k+1,l}^{m,c}$ the complement of
the set $G_{k+1,l}^{m}.$ It is easy to see that, for $0\leq t\leq \lambda ,$ 
\begin{equation}
E\exp \left\{ t\left| \widetilde{S}_{k,j}^{m}\right| \right\} =Q_{1}+Q_{2},
\label{P.1}
\end{equation}
where 
\begin{eqnarray}
Q_{1} &=&E\exp \left\{ t\left| \widetilde{S}_{k,j}^{m}\right| \right\} 
\mathbf{1}\left( G_{k+1,2j-1}^{m,c}\cup G_{k+1,2j-1}^{m,c}\right) ,
\label{P.2} \\
Q_{2} &=&E\exp \left\{ t\left| \widetilde{S}_{k,j}^{m}\right| \right\} 
\mathbf{1}\left( G_{k+1,2j-1}^{m}\cap G_{k+1,2j-1}^{m}\right) {.}
\label{P.2-2-BB}
\end{eqnarray}

First we give an estimate for $Q_{1}.$ Applying H\"{o}lder's inequality, we
obtain from (\ref{P.2}), 
\begin{equation}
Q_{1}\leq \left( \exp \left\{ 2t\left| \widetilde{S}_{k,j}^{m}\right|
\right\} \right) ^{1/2}\left( P\left( G_{k+1,2j-1}^{m,c}\right)
^{1/2}+P\left( G_{k+1,2j}^{m,c}\right) ^{1/2}\right) .  \label{P.4}
\end{equation}
By Lemma \ref{L P1} we have with $l=2j-1,2j$ 
\begin{equation}
P\left( G_{k+1,l}^{m,c}\right) =P\left( \left| \widetilde{Y}%
_{k+1,l}^{m}\right| >c^{\prime \prime }B_{k+1,l}^{m}\right) \leq 2\exp
\left\{ -c_{2}B_{k+1,l}^{m}\right\} .  \label{P.5-1}
\end{equation}
Note that by Proposition \ref{P B2}, we have $c_{3}^{-1}\leq
B_{k+1,2j-1}^{m}/B_{k+1,2j}^{m}\leq c_{3},$ which implies $B_{k+1,l}\geq
c_{4}B_{k,j}$ for $l=2j-1,2j.$ Then from (\ref{P.5-1}) it follows that 
\begin{equation}
P\left( G_{k+1,l}^{m,c}\right) \leq 2\exp \left\{ -c_{5}B_{k,j}^{m}\right\}
,\quad l=2j-1,2j.  \label{P.5}
\end{equation}

Inserting the bound provided by Lemma \ref{Lemma-exp-bound} and the
inequality (\ref{P.5}) into (\ref{P.4}) and choosing $t$ sufficiently small
we obtain 
\begin{equation*}
Q_{1}\leq c_{6}\exp \left\{ \left( c_{7}t^{2}-c_{8}\right)
B_{k,j}^{m}\right\} \leq c_{6}\exp \left\{ -\frac{1}{2}c_{8}B_{k,j}^{m}%
\right\} \leq c_{6}.
\end{equation*}

Now we shall give a bound for $Q_{2}.$ Recall that the r.v.'s $\widetilde{Y}%
_{k+1,l}^{m},$ $l=2j-1,2j$ are smooth (belong to the class $\mathfrak{D}(r)$%
), by Lemma \ref{L B3}. By virtue of Lemma \ref{L B5} and of the assumption $%
B_{k+1,2j-1}^{m}\geq c^{\prime }$ and $B_{k+1,2j}^{m}\geq c^{\prime },$ on
the set $G_{k+1,2j-1}^{m}\cap G_{k+1,2j}^{m}$ we have 
\begin{equation}
|\widetilde{S}_{k,j}^{m}|\leq c_{9}\left\{
1+U_{k+1,2j-1}^{m}+U_{k+1,2j}^{m}\right\} ,  \label{P.8b}
\end{equation}
where for $l=2j-1,2j$ 
\begin{equation*}
U_{k+1,l}^{m}=(\widetilde{Y}_{k+1,l}^{m,\ast })^{2}/B_{k+1,l}^{m},\quad 
\widetilde{Y}_{k+1,l}^{m,\ast }=\widetilde{Y}_{k+1,l}^{m}\mathbf{1}\left(
\left| \widetilde{Y}_{k+1,l}^{m}\right| \leq B_{k+1,l}^{m}\right) .
\end{equation*}
According to (\ref{P.2-2-BB}) and (\ref{P.8b}) 
\begin{eqnarray}
Q_{2} &\leq &E\exp \left\{ tc_{10}\left(
1+U_{k+1,2j-1}^{m}+U_{k+1,2j}^{m}\right) \right\}   \notag \\
&=&\exp \left\{ tc_{10}\right\} E\exp \left\{
tc_{10}U_{k+1,2j-1}^{m}\right\} E\exp \left\{ tc_{6}U_{k+1,2j}^{m}\right\} .
\label{P.9}
\end{eqnarray}
By Lemma \ref{L A3} (see the Appendix) we have 
\begin{equation}
E\exp \left\{ c_{10}U_{k+1,2j-1}^{m}\right\} \leq 1+2/c_{10}  \label{P.9aa}
\end{equation}
and a similar bound holds true for $U_{k+1,2j-1}^{m}.$ Taking $t$
sufficiently small, from (\ref{P.9}) and (\ref{P.9aa}) we obtain $Q_{2}\leq
c_{11}.$ Combining the estimates $Q_{1}\leq c_{6}$ and $Q_{2}\leq c_{11}$
obtained above with (\ref{P.1}) yields the lemma.
\end{proof}

\begin{lemma}
\label{L P3}For any $m=0,...,M$%
\begin{equation*}
E\exp \left\{ c_{1}\left| \widetilde{S}_{0}^{m}\right| \right\} \leq c_{2}.
\end{equation*}
\end{lemma}

\begin{proof}
The argument is similar to that for Lemma \ref{L P2}, and therefore will not
be given here. The only difference is that instead of Lemma \ref{L B5} we
make use of Lemma \ref{L B4}.
\end{proof}

Now Lemma \ref{lemma-basic} follows easily from Lemmas \ref{L P2}, \ref{L P3}
and Lemma \ref{L A1} in the Appendix.

\subsection{Proof of Theorem \ref{T3-2}}

The idea of the proof is to decompose the function $f$ into a Haar expansion
and then to make use of the closeness properties of the sequences $%
\widetilde{X}_{i},$ $i=1,...,n$ and $N_{i},$ $i=1,...,n$ over the dyadic
blocks. For this the representation provided by Proposition \ref{P R1} and
the exponential inequalities in Lemma \ref{lemma-basic} are crucial.

For the sake of brevity set 
\begin{equation*}
S_{n}(f)=\sum_{i=1}^{n}f(t_{i})\left( \widetilde{X}_{i}-N_{i}\right) .
\end{equation*}
What we have to show is that for any $t$ satisfying $\left| t\right| \leq
c_{0},$%
\begin{equation}
E\exp \{t(\log n)^{-2}S_{n}(f\})\leq \exp \left\{ t^{2}c_{1}\right\} .
\label{Z.1}
\end{equation}
Toward this end let $M=[\log _{2}(n/n_{0})]$ and note that according to
Proposition \ref{P R1}, 
\begin{equation*}
S_{n}(f)=\sum_{m=0}^{M}S^{m}\quad \text{where\quad }S^{m}=\sum_{i\in
J_{m}}f\left( t_{i}^{m}\right) \left( \widetilde{Y}_{i}^{m}-W_{i}^{m}\right)
.
\end{equation*}
By H\"{o}lder's inequality 
\begin{equation}
E\exp \left\{ t(\log n)^{-2}S_{n}(f)\right\} \leq \prod_{m=0}^{M}\left(
E\exp \left\{ t(M+1)(\log n)^{-2}S^{m}\right\} \right) ^{1/(M+1)}.
\label{Z.4}
\end{equation}
Set for brevity 
\begin{equation}
u_{n}=(M+1)(\log n)^{-2}.  \label{Z.5}
\end{equation}
Obviously $u_{n}\leq 1$ for $n$ large enough (such that $\log n\geq 2$).

It is easy to see that inequality (\ref{Z.1}) follows from (\ref{Z.4}) if we
prove that for $m=0,...,M$ and any $t$ satisfying $\left| t\right| \leq
c_{0} $%
\begin{equation}
E\exp \left\{ tu_{n}S^{m}\right\} \leq \exp \left\{ t^{2}c_{1}\right\} .
\label{Z.6}
\end{equation}
In the sequel we will give a proof of (\ref{Z.6}).

First we consider the case $m=0$. By H\"{o}lder's inequality, 
\begin{equation}
E\exp \left\{ tu_{n}S^{0}\right\} \leq \left( E\exp \left\{
2tu_{n}\sum_{i\in J_{0}}f(t_{i}^{0})\widetilde{Y}_{i}^{0}\right\} E\exp
\left\{ 2tu_{n}\sum_{i\in J_{0}}f(t_{i}^{0})W_{i}^{0}\right\} \right) ^{1/2}.
\label{Z.8}
\end{equation}
Since $\widetilde{Y}_{i}^{0},$ $i\in J_{0}$ are independent we have 
\begin{equation*}
E\exp \left\{ 2tu_{n}\sum_{i\in J_{0}}f(t_{i}^{0})\widetilde{Y}%
_{i}^{0}\right\} =\prod_{i\in J_{0}}E\exp \left\{ 2tu_{n}f(t_{i}^{0})%
\widetilde{Y}_{i}^{0}\right\} .
\end{equation*}
By choosing the constant $c_{0}$ small enough we can easily guarantee that $%
\left| 2tu_{n}f(t_{i}^{0})\right| \leq \lambda /3$, and by Lemma \ref{L A1}
we obtain 
\begin{equation}
E\exp \left\{ 2tu_{n}\sum_{i\in J_{0}}f(t_{i}^{0})\widetilde{Y}%
_{i}^{0}\right\} \leq \exp \left\{ c_{2}t^{2}\sum_{i\in J_{0}}E(\widetilde{Y}%
_{i}^{0})^{2}\right\} .  \label{Z.7a}
\end{equation}
Since $E(\widetilde{Y}_{i}^{0})^{2}=E^{\prime }(X_{i2^{M}})^{2}\leq C_{\max
} $ for $i\in J_{0},$ and $\#J_{0}\leq 2n_{\min }$ (see Section \ref
{sec:SNot}), we have $\sum_{i\in J_{0}}E(\widetilde{Y}_{i}^{0})^{2}\leq
c_{3},$ which in conjunction with (\ref{Z.7a}) yields 
\begin{equation*}
E\exp \{2tu_{n}\sum_{i\in J_{0}}f(t_{i}^{0})\widetilde{Y}_{i}^{0}\}\leq \exp
\left\{ c_{4}t^{2}\right\} .
\end{equation*}
An analogous bound holds true for the second expectation in (\ref{Z.8}).
From these bounds and from (\ref{Z.8}) we obtain (\ref{Z.6}) for $m=0.$

For the case $m\geq 1$ introduce the function $g(s)=f(a(s)),$ $s\in \lbrack
0,1],$ where $a(s)$ is defined by (\ref{BCN.4a}). Set for brevity $%
s_{i}^{m}=b(t_{i}^{m}),$ $i\in J_{m}$. Then for the sum $S^{m}$ we get the
following representation: 
\begin{equation*}
S^{m}=\sum_{i\in J_{m}}g\left( s_{i}^{m}\right) \left( \widetilde{Y}%
_{i}^{m}-W_{i}^{m}\right) .
\end{equation*}
Let $g_{m}$ be the truncated Haar expansion of $g$ for $m\geq 1$ (see (\ref
{HE.5}): 
\begin{equation}
g_{m}=c_{0}(g)h_{0}+\sum_{k=0}^{m-1}\sum_{j=1}^{2^{k}}c_{k,j}(g)h_{k,j},
\label{Y.2}
\end{equation}
where $c_{0}(g)$ and $c_{k,j}(g)$ are the corresponding Fourier-Haar
coefficients defined by (\ref{HE.4}) with $g$ replacing $f$. Then obviously 
\begin{equation*}
S^{m}=S_{1}^{m}+S_{2}^{m}
\end{equation*}
where 
\begin{eqnarray}
S_{1}^{m} &=&\sum_{i\in J_{m}}\left( g\left( s_{i}^{m}\right) -g_{m}\left(
s_{i}^{m}\right) \right) \left( \widetilde{Y}_{i}^{m}-W_{i}^{m}\right) ,
\label{Y.3a} \\
S_{2}^{m} &=&\sum_{i\in J_{m}}g_{m}\left( s_{i}^{m}\right) \left( \widetilde{%
Y}_{i}^{m}-W_{i}^{m}\right) .  \notag
\end{eqnarray}
By H\"{o}lder's inequality 
\begin{equation}
E\exp \left\{ tu_{n}S^{m}\right\} \leq \left( E\exp \left\{
2tu_{n}S_{1}^{m}\right\} E\exp \left\{ 2tu_{n}S_{2}^{m}\right\} \right)
^{1/2}.  \label{Y.4}
\end{equation}
Now the inequality (\ref{Z.6}) for $m\geq 1$ will be established if we prove
that both expectations on the right-hand side of (\ref{Y.4}) are bounded by $%
\exp \left\{ t^{2}c\right\} $. These inequalities are the subject of
Propositions \ref{P M1} and \ref{P M2} below. This completes the proof of
Theorem \ref{T3-2}.

First we prove the bound for the first expectation on the right hand side of
(\ref{Y.4}).

\begin{proposition}
\label{P M1}For any $m=1,...,M$ and $t$ satisfying $\left| t\right| \leq
c_{0}$ we have 
\begin{equation*}
E\exp \{tu_{n}S_{1}^{m}\}\leq \exp \left\{ t^{2}c_{1}\right\} .
\end{equation*}
\end{proposition}

\begin{proof}
Since by (\ref{BCN.5b}) the function $a(s)$ is Lipschitz and $f\in \mathcal{H%
}(\frac{1}{2},L),$ it is easy to see that the function $g(s)=f(a(s))$ is
also in a H\"{o}lder ball $\mathcal{H}(\frac{1}{2},L_{0})$ but with another
absolute constant $L_{0}.$ By H\"{o}lder's inequality 
\begin{equation}
E\exp \{tu_{n}S_{1}^{m}\}\leq \left( E\exp \left\{ \sum_{i\in J_{m}}\rho _{i}%
\widetilde{Y}_{i}^{m}\right\} E\exp \left\{ -\sum_{i\in J_{m}}\rho
_{i}W_{i}^{m}\right\} \right) ^{1/2},  \label{Y.5}
\end{equation}
where $\rho _{i}=2tu_{n}(g(s_{i}^{m})-g_{m}(s_{i}^{m}))$ and $\left|
t\right| \leq c_{0}$ for some sufficiently small absolute constant $c_{0}.$
Note that by Proposition \ref{P HE2} we have $\left\| g-g_{m}\right\|
_{\infty }\leq L_{0}2^{-m/2}.$ Therefore for $\left| t\right| \leq c_{0}$
(where $c_{0}$ is small) 
\begin{equation*}
\left| \rho _{i}\right| \leq c_{1}\left| t\right| u_{n}2^{-m/2}\leq
c_{1}\left| t\right| 2^{-m/2}\leq \lambda /3.
\end{equation*}
Then according to Lemma \ref{L A1} we get for $i\in J_{m}$%
\begin{equation}
E\exp \left\{ \rho _{i}\widetilde{Y}_{i}^{m}\right\} \leq \exp \left\{ \rho
_{i}^{2}E(\widetilde{Y}_{i}^{m})^{2}\right\} \leq \exp \left\{
c_{2}t^{2}2^{-m}E(X_{i}^{m})^{2}\right\} .  \label{Y.7a}
\end{equation}
An analogous bound holds true for the normal r.v.'s $W_{i}^{m},$ $i\in J_{m}:
$%
\begin{equation}
E\exp \left\{ -\rho _{i}W_{i}^{m}\right\} \leq \exp \left\{
c_{2}t^{2}2^{-m}E(X_{i}^{m})^{2}\right\} .  \label{Y.7b}
\end{equation}
Taking into account that $\widetilde{Y}_{i}^{m},$ $i\in J_{m}$ and $%
W_{i}^{m},$ $i\in J_{m}$ are sequences of independent r.v.'s and inserting (%
\ref{Y.7a}) and (\ref{Y.7b}) into (\ref{Y.5}), we obtain 
\begin{equation}
E\exp \{tu_{n}S_{1}^{m}\}\leq \exp \left\{ c_{3}t^{2}2^{-m}\sum_{i\in
J_{m}}E(X_{i}^{m})^{2}\right\} .  \label{Y.8}
\end{equation}
Now we remark that $\#J_{m}\leq 2^{m+1}.$ Hence by (\ref{NR.1}) 
\begin{equation}
\sum_{i\in J_{m}}E(X_{i}^{m})^{2}\leq \#J_{m}C_{\max }\leq 2^{m+1}C_{\max }.
\label{Y.9}
\end{equation}
Inserting (\ref{Y.9}) into (\ref{Y.8}), we obtain the result.
\end{proof}

Now we will find the bound for the second expectation on the right hand side
of (\ref{Y.4}).

\begin{proposition}
\label{P M2}For any $m=1,...,M$ and $t$ satisfying $\left| t\right| \leq
c_{0}$ we have 
\begin{equation*}
E\exp \{tu_{n}S_{2}^{m}\}\leq \exp \left\{ t^{2}c_{1}\right\} .
\end{equation*}
\end{proposition}

\begin{proof}
From (\ref{Y.3a}), (\ref{Y.2}) and (\ref{HE.2x-1}) we obtain 
\begin{equation*}
S_{2}^{m}=c_{0}(g)\left( \widetilde{Y}_{0,1}^{m}-W_{0,1}^{m}\right)
+\sum_{k=0}^{m-1}2^{k/2}\sum_{j=1}^{2^{k}}c_{k,j}(g)\left( T_{k,j}^{\ast
,m}-V_{k,j}^{\ast ,m}\right) 
\end{equation*}
where 
\begin{equation}
T_{k,j}^{\ast ,m}=\widetilde{Y}_{k+1,2j-1}^{m}-\widetilde{Y}%
_{k+1,2j}^{m},\quad V_{k,j}^{\ast ,m}=W_{k+1,2j-1}^{m}-W_{k+1,2j}^{m}
\label{X.8-1}
\end{equation}
(compare with (\ref{QI.1x-1}) and (\ref{QI.1x-2})). Here $\widetilde{Y}%
_{k,j}^{m}$ and $W_{k,j}^{m}$ are defined by (\ref{YYY}). Set in analogy to (%
\ref{Aux-st-1}) 
\begin{equation}
S_{0}^{m}=\widetilde{Y}_{0,1}^{m}-W_{0,1}^{m},\quad
S_{k,j}^{m}=T_{k,j}^{\ast ,m}-V_{k,j}^{\ast
,m},\;j=1,...,2^{k},\;k=0,...,m-1.  \label{X.8}
\end{equation}
Since the function $g(s)$ is in the H\"{o}lder ball with a H\"{o}lder
constant $L_{0},$ according to Proposition \ref{P HE1} we have the following
bounds for the Fourier-Haar coefficients: 
\begin{equation}
c_{0}(g)\leq L_{0}/2,\quad \left| c_{k,j}(g)\right| \leq
2^{-3/2}L_{0}2^{-k},\quad j=1,...,2^{k},\;k=0,...,m-1.  \label{X.2}
\end{equation}
Note also that by Lemma \ref{lemma-basic} there is an absolute constant $%
t_{0}$ sufficiently small such that for $|v|\leq t_{0}$ 
\begin{equation}
E\exp \left\{ v\widetilde{S}_{0}^{m}\right\} \leq \exp \left\{
c_{1}v^{2}\right\} ,\quad E\exp \left\{ v\widetilde{S}_{k,j}^{m}\right\}
\leq \exp \left\{ c_{1}v^{2}\right\}   \label{X.3}
\end{equation}
for $j=1,...,2^{k}$ and $k=0,...,m-1,$ where $\widetilde{S}_{0}^{m}$ and $%
\widetilde{S}_{k,j}^{m}$ are defined by (\ref{Aux-st-1}).

By H\"{o}lder's inequality we have, for any $t$ satisfying $\left| t\right|
\leq c_{0}\leq t_{0},$%
\begin{equation}
E\exp \left\{ tu_{n}S_{2}^{m}\right\} \leq \left( E\exp \left\{
t(m+1)u_{n}c_{0}(g)S_{0}^{m}\right\} \prod_{k=0}^{m-1}E\exp \left\{
t(m+1)u_{n}U_{k}\right\} \right) ^{1/(m+1)},  \label{XX.0}
\end{equation}
where 
\begin{equation}
U_{k}=2^{k/2}\sum_{j=1}^{2^{k}}c_{k,j}(g)S_{k,j}^{m},\quad k=0,...,m-1.
\label{X.4}
\end{equation}
The claim will be established, if we show that the constant $c_{0}$ can be
chosen such that for $t$ satisfying $\left| t\right| \leq c_{0},$ 
\begin{equation}
E\exp \left\{ tu_{m,n}c_{0}(g)S_{0}^{m}\right\} \leq \exp \left\{
c_{2}t^{2}\right\}   \label{X.5}
\end{equation}
and 
\begin{equation}
E\exp \left\{ tu_{m,n}U_{k}\right\} \leq \exp \left\{ c_{2}t^{2}\right\} ,
\label{X.6}
\end{equation}
where for the sake of brevity we set $u_{m,n}=(m+1)u_{n}.$

It is easy to show (\ref{X.5}). For this we note that by (\ref{X.2}) and (%
\ref{Z.5}), for $\left| t\right| \leq c_{0}$ we have 
\begin{equation}
\left| tu_{m,n}c_{0}(g)\right| \leq c_{3}\left| t\right|
(m+1)(M+1)L_{0}/\log ^{2}n\leq c_{4}c_{0}\leq t_{0},  \label{X.7}
\end{equation}
if the constant $c_{0}$ is small enough. Then the inequality (\ref{X.5})
follows from (\ref{X.3}) and from (\ref{X.7}).

The proof of (\ref{X.6}) is somewhat more involved. The main problem is that 
$S_{k,j}^{m},$ $j=1,...,2^{k}$ are dependent and therefore we cannot make
use of the product structure of the exponent $\exp \left\{ tU_{k}\right\} $
directly. However Proposition \ref{P B1} ensures that the components of the
sum $U_{k}$ (see (\ref{X.4})) are \textit{almost }independent, which allows
to exploit the product structure in an implicit way. The main idea is to
''substitute'' $S_{k,j}^{m},$ $j=1,...,2^{k}$ by $\widetilde{S}_{k,j}^{m},$ $%
j=1,....,2^{k}$ which are independent. With this in mind we write 
\begin{equation*}
U_{k}=U_{k}^{1}+U_{k}^{2},
\end{equation*}
where 
\begin{equation*}
U_{k}^{1}=2^{k/2}\sum_{j=1}^{2^{k}}c_{k,j}(g)\widetilde{S}_{k,j}^{m},\quad
U_{k}^{2}=2^{k/2}\sum_{j=1}^{2^{k}}c_{k,j}(g)\left( S_{k,j}^{m}-\widetilde{S}%
_{k,j}^{m}\right) .
\end{equation*}
Then by H\"{o}lder's inequality, 
\begin{equation}
E\exp \left\{ tu_{m,n}U_{k}\right\} \leq \left( E\exp \left\{
2tu_{m,n}U_{k}^{1}\right\} E\exp \left\{ 2tu_{m,n}U_{k}^{2}\right\} \right)
^{1/2}.  \label{X.16}
\end{equation}
Now we proceed to estimate the first expectation on the right-hand side of (%
\ref{X.16}). We make use of the independence of $\widetilde{S}_{k,j}^{m},$ $%
j=1,...,2^{k}$ (see Lemma \ref{L B2}) to get 
%%%%%%%%%%%%%%%%%%%%%%%%%%%%%%%%%%
\begin{equation}
E\exp \left\{ 2tu_{m,n}U_{k}^{1}\right\} =\prod_{j=1}^{2^{k}}E\exp \left\{
tq_{j}\widetilde{S}_{k,j}^{m}\right\} ,  \label{X.17}
\end{equation}
where $q_{j}=q_{m,n,k,j}=2u_{m,n}2^{k/2}c_{k,j}(g).$ Note that by (\ref{X.2}%
) and (\ref{Z.5}) 
\begin{equation*}
\left| tq_{j}\right| \leq \left| 2tu_{m,n}2^{k/2}c_{k,j}(g)\right| \leq
c_{5}\left| t\right| 2^{-k/2}\leq t_{0},
\end{equation*}
provided $c_{0}$ is small enough. It then follows from (\ref{X.3}) that for $%
j=1,...,2^{k}$%
\begin{equation}
E\exp \left\{ tq_{j}\widetilde{S}_{k,j}^{m}\right\} \leq \exp \left\{
c_{6}t^{2}2^{-k}\right\} .  \label{X.18}
\end{equation}
Inserting (\ref{X.18}) into (\ref{X.17}) we find the bound 
\begin{equation}
E\exp \left\{ 2tu_{m,n}U_{k}^{1}\right\} \leq \exp \left\{
c_{7}t^{2}\right\} .  \label{XX.5}
\end{equation}
Thus we have estimated the first expectation on the right hand side of (\ref
{X.16}). It remains to estimate the second one.

Note that 
\begin{equation*}
\widetilde{S}_{k,j}^{m}-S_{k,j}^{m}=\left( \widetilde{T}_{k,j}^{m}-T_{k,j}^{%
\ast ,m}\right) -\left( V_{k,j}^{m}-V_{k,j}^{\ast ,m}\right) .
\end{equation*}
Hence 
\begin{equation*}
U_{k}^{2}=U_{k}^{3}+U_{k}^{4}
\end{equation*}
where 
\begin{equation*}
U_{k}^{3}=2^{k/2}\sum_{j=1}^{2^{k}}c_{k,j}(g)\left( \widetilde{T}%
_{k,j}^{m}-T_{k,j}^{\ast ,m}\right) ,\quad
U_{k}^{4}=2^{k/2}\sum_{j=1}^{2^{k}}c_{k,j}(g)\left(
V_{k,j}^{m}-V_{k,j}^{\ast ,m}\right) .
\end{equation*}
By H\"{o}lder's inequality we obtain 
\begin{equation}
E\exp \left\{ 2tu_{m,n}U_{k}^{2}\right\} \leq \left( E\exp \left\{
4tu_{m,n}U_{k}^{3}\right\} E\exp \left\{ 4tu_{m,n}U_{k}^{4}\right\} \right)
^{1/2}.  \label{XX.1}
\end{equation}
Since $\widetilde{T}_{k,j}^{m}-T_{k,j}^{\ast ,m},$ $j=1,...,2^{k}$ is a
sequence of independent r.v.'s, we get 
\begin{equation}
E\exp \left\{ 4tu_{m,n}U_{k}^{3}\right\} \leq \prod_{j=1}^{2^{k}}E\exp
\left\{ 2tq_{j}\left( \widetilde{T}_{k,j}^{m}-T_{k,j}^{\ast ,m}\right)
\right\}   \label{XX.2}
\end{equation}
where $q_{j}$ is defined above (see (\ref{X.17})). The definitions of $%
\widetilde{T}_{k,j}^{m}$ and of $T_{k,j}^{\ast ,m}$ (see (\ref{QI.1x-1}) and
(\ref{X.8-1})) imply 
\begin{equation*}
\widetilde{T}_{k,j}^{m}-T_{k,j}^{\ast ,m}=\beta _{2j}\widetilde{Y}%
_{k+1,2j-1}^{m}-\beta _{2j-1}\widetilde{Y}_{k+1,2j}^{m}.
\end{equation*}
Hereafter we abbreviate $\beta _{i}=\alpha _{k+1,i}^{m}-1,$ $%
B_{i}=B_{k+1,i}^{m}.$ Then 
\begin{equation}
E\exp \left\{ 2tq_{j}\left( \widetilde{T}_{k,j}^{m}-T_{k,j}^{\ast ,m}\right)
\right\} =E\exp \left\{ tq_{j}\beta _{2j}\widetilde{Y}_{k+1,2j-1}^{m}\right%
\} E\exp \left\{ -tq_{j}\beta _{2j-1}\widetilde{Y}_{k+1,2j}^{m}\right\} .
\label{XX.3}
\end{equation}
Since by Proposition \ref{P B2} $B_{2j}\leq c_{8}B_{2j-1},$ we have $\beta
_{2j}\leq 1+c_{8}.$ Hence by (\ref{X.2}) and (\ref{Z.5}) 
\begin{equation}
\left| tq_{j}\beta _{2j}\right| \leq c_{9}\left| t\right| 2^{-k/2}\beta
_{2j}\leq \lambda /3  \label{XX.6}
\end{equation}
for $t$ sufficiently small. By (\ref{YYY}) and by Lemma \ref{Lemma-CNSTR}, $%
\widetilde{Y}_{k+1,2j-1}^{m}$ is a sum of independent r.v.'s which satisfy
Sakhanenko's condition (\ref{NR.2}). Hence using Lemma \ref{L A1} we obtain 
\begin{eqnarray*}
E\exp \left\{ tq_{j}\beta _{2j}\widetilde{Y}_{k+1,2j-1}^{m}\right\} 
&=&\prod_{i\in I_{k+1,2j-1}^{m}}E\exp \left\{ tq_{j}\beta _{2j}\widetilde{Y}%
_{i}^{m}\right\}  \\
&\leq &\prod_{i\in I_{k+1,2j-1}^{m}}\exp \left\{ t^{2}q_{j}^{2}\beta
_{2j}^{2}E\left( \widetilde{Y}_{i}^{m}\right) ^{2}\right\} .
\end{eqnarray*}
By (\ref{XX.6}) 
\begin{eqnarray*}
E\exp \left\{ tq_{j}\beta _{2j}\widetilde{Y}_{k+1,2j-1}^{m}\right\}  &\leq
&\prod_{i\in I_{k+1,2j-1}^{m}}\exp \left\{ c_{10}t^{2}2^{-k/2}\beta
_{2j}^{2}E\left( \widetilde{Y}_{i}^{m}\right) ^{2}\right\}  \\
&=&\exp \left\{ c_{10}t^{2}2^{-k/2}\beta _{2j}^{2}B_{2j-1}\right\} .
\end{eqnarray*}
Taking into account Proposition \ref{P B1}, we obtain 
\begin{equation*}
\beta _{2j}^{2}B_{2j-1}=\left( \sqrt{B_{2j}}-\sqrt{B_{2j-1}}\right) ^{2}\leq
\left| B_{2j}-B_{2j-1}\right| \leq c_{11}.
\end{equation*}
This proves that 
\begin{equation*}
E\exp \left\{ tq_{j}\beta _{2j}\widetilde{Y}_{k+1,2j-1}^{m}\right\} \leq
\exp \left\{ c_{12}t^{2}2^{-k/2}\right\} .
\end{equation*}
For the second expectation on the right hand side of (\ref{XX.3}) we prove
an analogous bound. Invoking these bounds in (\ref{XX.3}) we get 
\begin{equation}
E\exp \left\{ 2tq_{j}\left( \widetilde{T}_{k,j}^{m}-T_{k,j}^{\ast ,m}\right)
\right\} \leq \exp \left\{ c_{13}t^{2}2^{-k/2}\right\} .  \label{XX.9}
\end{equation}
Inserting in turn (\ref{XX.9}) into (\ref{XX.2}) we arrive at 
\begin{equation*}
E\exp \left\{ 4tu_{m,n}U_{k}^{3}\right\} \leq \exp \left\{
c_{14}t^{2}\right\} .
\end{equation*}
In the same way we prove an inequality for $U_{k}^{4}.$ Then by (\ref{XX.1})
we have 
\begin{equation}
E\exp \left\{ 2tu_{m,n}U_{k}^{2}\right\} \leq \exp \left\{
c_{14}t^{2}\right\} .  \label{XX.10}
\end{equation}
From (\ref{X.16}), (\ref{XX.10}) and (\ref{XX.5}) we obtain inequality (\ref
{X.6}), this completing the proof of the proposition.
\end{proof}

\section{Appendix}

In the course of the reasoning we made use of the following simple auxiliary
results.

\begin{lemma}
\label{L A1}Let $\xi $ be a real valued r.v.\ with mean $0$ and finite
variance: $E\xi =0,$ $0<E\xi ^{2}<\infty .$ Assume that Sakhanenko's
condition 
\begin{equation*}
\lambda E\left| \xi \right| ^{3}\exp \{\lambda |\xi |\}\leq E\xi ^{2}
\end{equation*}
holds true for some $\lambda >0.$ Then for all $|t|\leq \lambda /3$ 
\begin{equation*}
E\exp \{t\xi \}\leq \exp \left\{ t^{2}E\xi ^{2}\right\} .
\end{equation*}
\end{lemma}

\begin{proof}
Let $\mu (t)=E\exp (t\xi )$ and $\psi (t)=\log \mu (t)$ be the moment and
cumulant generating functions respectively. The conditions of the lemma
imply that $\mu (t)\leq c_{1}$ for any real $|t|\leq \lambda /3.$ Using a
three term Taylor expansion we obtain for $0\leq \nu \leq 1$%
\begin{equation*}
\psi (t)=\psi (0)+\psi ^{\prime }(0)t+\psi ^{\prime \prime }(0)\frac{t^{2}}{2%
}+\psi ^{\prime \prime }{}^{\prime }(\nu t)\frac{t^{3}}{6}.
\end{equation*}
Note that $\psi (0)=0,$ $\psi ^{\prime }(0)=0,$ $\psi ^{\prime \prime
}(0)=E\xi ^{2}$ and $\mu (t)\geq 1$ by Jensen's inequality, while for the
third derivative we have for any real $s$ satisfying $|s|\leq \lambda /3,$ 
\begin{equation*}
\psi ^{\prime \prime \prime }(s)=\mu ^{\prime \prime \prime }(s)\mu
(s)^{-1}-3\mu ^{\prime \prime }(s)\mu ^{\prime }(s)\mu (s)^{-2}+2\mu
^{\prime }(s)^{3}\mu (s)^{-3}.
\end{equation*}
Using H\"{o}lder's inequality and $\mu (s)\geq 1$ we obtain the bound 
\begin{equation*}
\left| \psi ^{\prime \prime \prime }(s)\right| \leq 6E|\xi |^{3}\exp
(\lambda |\xi |).
\end{equation*}
Since $|t|\leq \lambda /3,$ by Sakhanenko's condition we have 
\begin{equation*}
0\leq \psi (t)\leq \frac{t^{2}}{2}E\xi ^{2}+t^{3}E|\xi |^{3}\exp (\lambda
|\xi |)\leq t^{2}E\xi ^{2}.\;
\end{equation*}
\end{proof}

\begin{lemma}
\label{L A2}Let $\xi $ be a real valued r.v.\ such that $E\xi =0$ and 
\begin{equation*}
E\exp \{\lambda |\xi |\}\leq c_{1}
\end{equation*}
for some $\lambda \geq 0$ and $c_{1}\geq 1$. Then for all $|t|\leq \lambda /2
$ we have 
\begin{equation*}
E\exp \{t\xi \}\leq \exp \{c_{2}t^{2}\},
\end{equation*}
where $c_{2}=4c_{1}/\lambda ^{2}$.
\end{lemma}

\begin{proof}
The argument is similar to Lemma \ref{L A1}. We use the same notations. A
two term Taylor expansion yields, for $0\leq \nu \leq 1,$%
\begin{equation*}
\psi (t)=\psi (0)+\psi ^{\prime }(0)t+\psi ^{\prime \prime }(\nu t)\frac{%
t^{2}}{2}.
\end{equation*}
Since $x^{2}\leq 2\exp (|x|)$ for any real $x,$ we have for any $s$
satisfying $|s|\leq \lambda /2$ 
\begin{eqnarray}
0 &\leq &\psi ^{\prime \prime }(s)=\mu (s)^{-2}\{E\xi ^{2}\exp (s\xi )-(E\xi
\exp (s\xi ))^{2}\}  \notag \\
\  &\leq &E\xi ^{2}\exp (s\xi )\leq E\xi ^{2}\exp (\frac{\lambda }{2}|\xi
|)\leq 8\frac{c_{1}}{\lambda ^{2}}.  \notag
\end{eqnarray}
Consequently 
\begin{equation*}
0\leq \psi (t)=\psi ^{\prime \prime }(\nu t)\frac{t^{2}}{2}\leq 4\frac{c_{1}%
}{\lambda ^{2}}t^{2}.
\end{equation*}
\end{proof}

\begin{lemma}
\label{L A3}Let $\xi _{i},$ $i=1,...,n$ be a sequence of independent r.v.'s
such that for all $i=1,...,n$ we have $E\xi _{i}=0,$ $0<E\xi _{i}^{2}<\infty 
$ and 
\begin{equation*}
\lambda E|\xi _{i}|^{3}\exp \left\{ \lambda |\xi _{i}|\right\} \leq E\xi
_{i}^{2}
\end{equation*}
for some positive constant $\lambda .$ Set $S_{n}=\xi _{1}+...+\xi _{n},$ $%
B_{n}^{2}=ES_{n}^{2}$ and $S_{n}^{\ast }=S_{n}\mathbf{1}\left( \left|
S_{n}\right| \leq B_{n}^{2}\right) .$ Then 
\begin{equation*}
E\exp \left\{ c_{1}(S_{n}^{\ast }/B_{n})^{2}\right\} \leq 1+2/c_{1},
\end{equation*}
where $c_{1}=\frac{1}{4}\min \left\{ \lambda /3,1/2\right\} .$
\end{lemma}

\begin{proof}
Denote 
\begin{equation*}
F(x)=P\left( (S_{n}^{\ast }/B_{n})^{2}>x\right) .
\end{equation*}
First we shall prove that 
\begin{equation}
F(x)\leq 2\exp \{-c_{2}x\},\quad x\geq 0,  \label{AS.2}
\end{equation}
where $c_{2}=2c_{1}.$ For this we note that 
\begin{equation*}
F(x)=P\left( S_{n}^{\ast }/B_{n}>\sqrt{x}\right) +P\left( S_{n}^{\ast
}/B_{n}<-\sqrt{x}\right) .
\end{equation*}
It suffices to estimate only the first probability on the right hand side of
the above equality; the second can be treated in the same way. If $%
x>B_{n}^{2}$ then 
\begin{equation*}
P\left( S_{n}^{\ast }/B_{n}>\sqrt{x}\right) =0;
\end{equation*}
thus there is nothing to prove in this case. Let $x\leq B_{n}^{2}.$ Denoting 
$t=2c_{2}\sqrt{x},$ one obtains 
\begin{eqnarray}
P\left( S_{n}^{\ast }>\sqrt{x}\right)  &\leq &P\left( S_{n}>\sqrt{x}\right)
\leq \exp \left\{ -t\sqrt{x}\right\} E\exp \left\{ tS_{n}/B_{n}\right\}  
\notag \\
\  &=&\exp \left\{ -t\sqrt{x}\right\} \prod_{i=1}^{n}E\exp \left\{ t\xi
_{i}/B_{n}\right\} .  \label{AS.3}
\end{eqnarray}
Note that $t/B_{n}=2c_{2}\sqrt{x}/B_{n}\leq 2c_{2}\leq \lambda /3.$ Hence by
Lemma \ref{L A1} 
\begin{equation*}
E\exp \left\{ t\xi _{i}/B_{n}\right\} \leq \exp \left\{ t^{2}E\xi
_{i}^{2}/B_{n}^{2}\right\} .
\end{equation*}
Inserting this into (\ref{AS.3}) we get 
\begin{eqnarray*}
P\left( S_{n}^{\ast }/B_{n}>\sqrt{x}\right)  &\leq &\exp \left\{ -t\sqrt{x}%
\right\} \prod_{i=1}^{n}\exp \left\{ t^{2}E\xi _{i}^{2}/B_{n}^{2}\right\}  \\
\  &=&\exp \left\{ -t\sqrt{x}+t^{2}\right\} \leq \exp \left\{
-c_{2}x\right\} 
\end{eqnarray*}
which proves (\ref{AS.2}). Integrating by parts we obtain 
\begin{eqnarray*}
E\exp \left\{ c_{1}(S_{n}^{\ast })^{2}/B_{n}\right\}  &=&\int_{0}^{\infty
}\exp \{c_{1}x\}dF(x) \\
\  &=&1+\int_{0}^{\infty }F(x)\exp \{c_{1}x\}dx \\
\  &\leq &1+2\int_{0}^{\infty }\exp \{c_{1}x-c_{2}x\}dx \\
\  &\leq &1+2/c_{1}.
\end{eqnarray*}
\end{proof}

\textbf{Acknowledgment.} We would like to thank A. Zaitsev for useful
discussions and the referee for remarks which helped to improve the
presentation of the paper.

%%%%%%%%%%%%%%%%%%%%%%%%%%%%%%%%%%%%%%%%%%%%%%%%%%%%%%%%%%%%%%%%%%%%%%%%%%%%%
%
%    References
%
%%%%%%%%%%%%%%%%%%%%%%%%%%%%%%%%%%%%%%%%%%%%%%%%%%%%%%%%%%%%%%%%%%%%%%%%%%%%%

\end{document}